\documentclass[12pt]{article}

\usepackage{amssymb}
\usepackage{latexsym}
\usepackage{amsmath}
\usepackage{amscd}
\usepackage{theorem}

\newtheorem{theorem}{Theorem}[section]
\newtheorem{proposition}[theorem]{Proposition}
\newtheorem{corollary}[theorem]{Corollary}
\newtheorem{lemma}[theorem]{Lemma}

{\theorembodyfont{\rmfamily}
\newtheorem{definition}[theorem]{Definition}
\newtheorem{remark}[theorem]{Remark}
\newtheorem{example}[theorem]{Example}

\newtheorem{choice}[theorem]{Choice}
}

\numberwithin{equation}{section}

\newenvironment{proof}
{\noindent {\em Proof.}}
{\hfill $\Box$}
\newcommand{\noProof} {\hfill$\Box$}

\newcommand{\eqdef}{\;{:=}\;}

\newcommand{\Q}{{\mathbb Q}}
\newcommand{\R}{{\mathbb R}}
\newcommand{\Z}{{\mathbb Z}}
\newcommand{\op}{\operatorname}
\newcommand{\mc}[1]{{\mathcal #1}}

\newcommand{\Ker}{\op{Ker}}

\newcommand{\Nov}{\op{Nov}}
\newcommand{\tensor}{\otimes}

\newcommand{\Crit}{\mc{C}}
\newcommand{\Critt}{\tilde{\Crit}}
\newcommand{\too}{\longrightarrow}

\newcommand{\Bas}{\mc{B}}
\newcommand{\Cnov}{CN}
\newcommand{\Eul}{\mc{E}}
\newcommand{\Tmorse}{T_m}
\newcommand{\Ttop}{T(\tilde{X})}

\title{Reidemeister torsion in generalized Morse theory}

\author{Michael Hutchings}

\date{}

\begin{document}

\maketitle

\begin{abstract}
In two previous papers with Yi-Jen Lee, we defined and computed a
notion of Reidemeister torsion for the Morse theory of closed 1-forms
on a finite dimensional manifold.  The present paper gives an a priori
proof that this Morse theory invariant is a topological invariant.  It
is hoped that this will provide a model for possible generalizations
to Floer theory.
\end{abstract}

In two papers with Yi-Jen Lee \cite{hl1,hl2}, we defined a notion of
Reidemeister torsion for the Morse theory of closed 1-forms on a
finite dimensional manifold. We consider the flow dual to the 1-form
via an auxiliary metric.  Our invariant, which we call $I$, multiplies
the algebraic Reidemeister torsion of the Novikov complex, which
counts flow lines between critical points, by a zeta function which
counts closed orbits of the flow.  For a closed 1-form in a rational
cohomology class, i.e. $d$ of a circle-valued function, we proved in
the above papers that $I$ equals a form of topological Reidemeister
torsion due to Turaev.  This implies {\em a posteriori} that $I$ is
invariant under homotopy of the circle-valued function and the
auxiliary metric.

In this paper we reprove these results using an opposite approach: we
first prove {\em a priori} that $I$ is a topological invariant,
depending only on the cohomology class of the closed 1-form.  We then
deduce that $I$ agrees with Turaev torsion, by using invariance to
reduce to the easier case of an exact 1-form.  This approach has two
advantages.  First, it works for closed 1-forms in an arbitrary
cohomology class, thus extending the results of our previous papers.
Second, and perhaps more importantly, the proof of invariance here
should provide a model for the possible construction of torsion in
Floer theory.

The contents of this paper are as follows.  In \S\ref{sec:intro} we
review the definition of the invariant $I$, state the main results,
and outline the proofs.  In \S\ref{sec:preparation} and
\S\ref{sec:bifurcationAnalysis} we prove that $I$ is invariant.  The
strategy is to study how $I$ varies in a generic one parameter family
of 1-forms and metrics.  In \S\ref{sec:preparation}, we prepare for
this analysis by classifying the bifurcations that generically occur,
and we also deal with the complication that infinitely many
bifurcations may occur in a finite time.  The heart of the paper is in
\S\ref{sec:bifurcationAnalysis}, where we analyze what happens in each
individual bifurcation.  While the torsion of the Novikov complex and
the zeta function can change, we will see that their product $I$ does
not.  In \S\ref{sec:comparison} we use topological invariance to give
a quick proof that $I$ agrees with Turaev torsion.  In
\S\ref{sec:conclusion} we discuss open questions and possible
generalizations.  Appendix~\ref{app:torsion} reviews algebraic aspects
of Reidemeister torsion that are used throughout the paper.
Appendix~\ref{app:euler} reviews how to remove a certain ambiguity in
Reidemeister torsion using Turaev's Euler structures.

\section{The invariant $I$}
\label{sec:intro}

We begin by reviewing the definition of the invariant $I$ from
\cite{hl2}.  We will emphasize geometric aspects which are important
for the present paper, and we generalize \cite{hl2} slightly by
allowing different covers in Choice~\ref{choice:covering}.  After
defining $I$, we will state the main results and outline the proofs.

\subsection{Setup and geometric definitions}
\label{sec:setup}

Let $X$ be a smooth, finite dimensional, closed (connected) manifold
with $\chi(X)=0$.  We consider a closed 1-form $\alpha$ and a
Riemannian metric $g$ on $X$.  Let $V\eqdef g^{-1}\alpha$ denote the
vector field dual to $\alpha$ via $g$.  We wish to count closed orbits
and flow lines of $V$, which are defined as follows.

A {\bf closed orbit} is a nonconstant map $\gamma:S^1\to X$ such that
$\gamma'(t)=-\lambda V_{\gamma(t)}$ for some $\lambda>0$.  There is a
minus sign because we work with the ``downward'' flow as in classical
Morse theory.  We consider two closed orbits to be equivalent if they
differ by a rotation of $S^1$.  The {\bf period} $p(\gamma)$ is the
largest integer $k$ such that $\gamma$ factors through a $k$-fold
covering $S^1\to S^1$.

For counting purposes, we can attach a sign to a generic closed orbit
as follows.  For $q\in\gamma(S^1)$, let $U\subset X$ be a hypersurface
intersecting $\gamma$ transversely near $q$, and let $\phi_q:U\to U$
be the return map (defined near $q$) which follows the flow $-V$ a
total of $p(\gamma)$ times around the image $\gamma(S^1)$.  The
linearized return map induces a map
\[
d\phi_q:T_qX/T_q\gamma(S^1)\to
T_qX/T_q\gamma(S^1)
\]
which does not depend on $U$.  The eigenvalues of
$d\phi_q$ do not depend on $q$.  We say that $\gamma$ is {\bf
nondegenerate} if $1-d\phi_q$ is invertible.  In this case we define the
{\bf Lefschetz sign} $(-1)^{\mu(\gamma)}$ to be the sign of
$\det(1-d\phi_q)$.

A {\bf critical point} is a zero of $\alpha$.  A critical point $p\in
X$ is {\bf nondegenerate} if the graph of $\alpha$ in the cotangent
bundle $T^*X$ intersects the zero section transversely at $p$.  In
this case the derivative $\nabla V:T_pX\to T_pX$ is invertible and
self-adjoint; the {\bf index} of $p$, denoted by $\op{ind}(p)$, is the
number of negative eigenvalues of $\nabla V$.  The {\bf descending
manifold} $\mc{D}(p)$ is the set of all $x\in X$ such that the
trajectory of the flow $+V$ starting at $x$ converges to
$p$. Similarly, the {\bf ascending manifold} $\mc{A}(p)$ is the set of
all $x\in X$ from which the trajectory of $-V$ converges to $p$.  If
$p$ is nondegenerate, then $\mc{D}(p)$ and $\mc{A}(p)$ are embedded
open balls of dimension $\op{ind}(p)$ and $n-\op{ind}(p)$,
respectively.

If $p$ and $q$ are critical points, a {\bf flow line} (of $-V$) from
$p$ to $q$ is a map $f:\R\to X$ such that $f'(t)=-V_{f(t)}$ and
$\lim_{t\to-\infty}f(t)=p$ and $\lim_{t\to+\infty}f(t)=q$.  We
consider two flow lines to be equivalent if they differ by a
translation of $\R$.  The space of flow lines from $p$ to $q$ is
naturally identified with $(\mc{D}(p)\cap\mc{A}(q))/\R$.  Thus, if $p$
and $q$ are nondegenerate, the expected dimension of the space of flow
lines from $p$ to $q$ is $\op{ind}(p)-\op{ind}(q)-1$.

We will impose the following transversality conditions.

\begin{definition}
\label{def:admissible}
The pair $(\alpha,g)$ is {\bf admissible} if:
\begin{description}
\item{(a)} All critical points of $V$ are nondegenerate.
\item{(b)} The ascending and descending manifolds of critical points
of $V$ intersect transversely. 
\item{(c)} All closed orbits of $V$ are nondegenerate.
\end{description}
\end{definition}
A straightforward transversality calculation
(cf. \cite{a-b,schwarz,thesis}) shows that for a fixed cohomology
class $[\alpha]\in H^1(X;\R)$, these conditions hold for generic pairs
$(\alpha,g)$.

\subsection{Coverings and Novikov rings}

To enable finite counting of closed orbits and critical points, we
work in a covering of $X$.

\begin{choice}
\label{choice:covering}
We choose a connected abelian covering $\pi:\tilde{X}\to X$ such that
$\pi^*\alpha$ is exact.
\end{choice}

We do not know how to remove the assumption that $\pi$ is abelian,
cf. \S\ref{sec:conclusion}.

Let $H$ denote the group of covering transformations.  There is a
surjection $H_1(X)\to H$, whose kernel consists of homology classes of
loops in $X$ that lift to $\tilde{X}$.

Our counting will take place in the {\bf Novikov ring}
$\Lambda\eqdef\op{Nov}(H;[-\alpha])$.  The meaning of this notation is
that if $G$ is an abelian group and $N:G\to\R$ is a homomorphism, then
$\op{Nov}(G;N)$ consists of possibly infinite formal sums $\sum_{g\in
G}a_g\cdot g$, with $a_g\in\Z$, such that for each $R\in\R$, there are
only finitely many nonzero coefficients $a_g$ with $N(g)<R$.  This
ring has the obvious addition, and the convolution product
\cite{novikov,hofer-salamon}.

\begin{example}
\label{ex:cyclic}
Suppose $\alpha=df$, where $f:X\to S^1$ is not nullhomotopic.  The
simplest choice is to take the covering $\tilde{X}$ to be a component
of the fiber product of $X$ and $\R$ over $S^1$.  Then $H\simeq\Z$,
and the Novikov ring is $\Lambda\simeq \Z((t))=\{\sum_{k=m}^\infty
a_kt^k|m,a_k\in\Z\}$, the ring of integer Laurent series.

This is essentially the setup of \cite{hl1}.  (In \cite{hl1},
$\tilde{X}$ was the entire fiber product of $X$ and $\R$ over $S^1$.
As a result, $t$ here equals $t^N$ in that paper, where $N$ is the
number of components of the fiber product, or equivalently the
divisibility of $[\alpha]\in H^1(X;\Z)$.)  For more refined
invariants, one can take $\tilde{X}$ to be the universal abelian
cover, as in \cite{hl2}.
\end{example}

\subsection{Counting closed orbits}
\label{sec:zeta}

If $(\alpha,g)$ is admissible, we count closed orbits using the {\bf
zeta function}
\begin{equation}
\label{eqn:zetaDef}
\zeta\eqdef\exp\sum_{\gamma\in\mc{O}}
\frac{(-1)^{\mu(\gamma)}}{p(\gamma)}[\gamma]\in\Lambda.
\end{equation}
(Cf. \cite{fried:homological,pajitnov:new}.)  Here $\mc{O}$ denotes
the set of closed orbits, and $[\gamma]\in H$ denotes the image of the
homology class $\gamma_*[S^1]$ under the projection $H_1(X)\to H$.

Let us review why $\zeta$ is a well defined element of $\Lambda$, as
the ideas in this argument will be important later.  First, we claim
that for each $R\in\R$, there are only finitely many closed orbits
$\gamma$ with $[-\alpha](\gamma)<R$.  Since $\alpha$ is closed, the
length of such an orbit away from the critical points is bounded above
by some multiple of $R$.  An elementary compactness argument
\cite{thesis,salamon:novikov} then shows that an infinite sequence of
such orbits would accumulate to either (i) a degenerate closed orbit,
or (ii) a ``broken'' closed orbit with stopovers at one or more
critical points.  Situation (i) would violate admissibility condition
(c).  In situation (ii), there would necessarily be a flow line from a
critical point of index $i$ to a critical point of index $\ge i$.
This would violate admissibility condition (b), since the expected
dimension of the space of such flow lines is negative.

Let $\Lambda^+$ denote the set of sums $\sum_{h\in H}a_h\cdot
h\in\Lambda$ such that $a_h=0$ whenever $[-\alpha](h)\le 0$.  Let
$\Lambda^+_\Q\eqdef\Lambda^+\tensor\Q$.  The above paragraph shows
that
\[
\sum_{\gamma\in\mc{O}}\frac{(-1)^{\mu(\gamma)}}{p(\gamma)}[\gamma] \in
\Lambda^+_\Q.
\]
Now $\exp:\Lambda^+_\Q\to 1+\Lambda^+_\Q$ is well defined by the usual
power series.   Therefore
$\zeta\in\Lambda\tensor\Q$.

But in fact $\zeta$ has integer coefficients, because there is a
product formula
\begin{equation}
\label{eqn:productFormula}
\zeta=\prod_{\gamma\in\mc{I}}(1\pm[\gamma])^{\pm1}.
\end{equation}
Here $\mc{I}$ denotes the set of {\bf irreducible} (period 1) closed
orbits, and the two signs associated to each irreducible orbit are
determined by the eigenvalues of the return map.  The formula is
proved by taking the logarithm of both sides,
cf. \cite{fried:survey,hl1,ionel-parker,smale}.

\begin{remark}
\label{remark:powerSeries}
The inverse of
$\exp$ above is also well defined via the usual power series
$\log(1+x)=\sum_{k=1}^\infty\frac{(-1)^{k+1}x^k}{k}$.
We will use this fact in \S\ref{sec:cyclic}.
\end{remark}

\subsection{Counting flow lines} 
\label{sec:morse}

We count flow lines using the {\bf Novikov complex}
$(\Cnov_*,\partial)$, which is defined as follows.  Let $\Critt_i$
denote the set of index $i$ critical points of $\pi^*V$ in
$\tilde{X}$.  Choose $f:\tilde{X}\to\R$ with $df=\pi^*\alpha$.  We
define $\Cnov_i$ to be the set of formal sums
$\sum_{\tilde{p}\in\Critt_i}a_{\tilde{p}}\cdot \tilde{p}$ with
$a_{\tilde{p}}\in\Z$, such that for each $R\in\R$, there are only
finitely many nonzero coefficients $a_{\tilde{p}}$ with
$f(\tilde{p})>R$.  The action of $H$ on $\Critt_i$ by covering
transformations makes $\Cnov_i$ into a module over the Novikov ring
$\Lambda$.  This module is free; one can obtain a basis by choosing a
lift of each critical point in $X$ to $\tilde{X}$.

We define the boundary operator $\partial:\Cnov_i\to\Cnov_{i-1}$ by
\begin{equation}
\partial\tilde{p}\eqdef
\sum_{\tilde{q}\in\Critt_{i-1}}
\langle\partial\tilde{p},\tilde{q}\rangle\cdot\tilde{q}
\end{equation}
for $\tilde{p}\in\Critt_i$.  Here
$\langle\partial\tilde{p},\tilde{q}\rangle$ denotes the
signed number of flow lines of $-\pi^*V$ from $\tilde{p}$ to
$\tilde{q}$.

The signs are determined as follows \cite{salamon:lectures}.  We
choose orientations of the descending manifolds of the critical points
in $X$, and lift them to orient the descending manifolds in
$\tilde{X}$.  Given a nondegenerate flow line from $\tilde{p}$ to
$\tilde{q}$, let $v\in T_{\tilde{p}}\mc{D}(\tilde{p})$ be an outward
tangent vector of the flow line.  The flow near the flow line
determines an isomorphism $T_{\tilde{p}}\mc{D}(\tilde{p})/v\to
T_{\tilde{q}}\mc{D}(\tilde{q})$.  We declare the flow line to have
positive sign if the orientations on $T_{\tilde{p}}\mc{D}(\tilde{p})$
and $\R v\oplus T_{\tilde{q}}\mc{D}(\tilde{q})$ agree.  (We do not
need to assume that $X$ is oriented.)

A compactness argument as in \S\ref{sec:zeta} and
\cite{salamon:novikov,pozniak,thesis}, using the fact that
$\pi^*\alpha$ is exact, shows that $\partial$ is well defined if
$(\alpha,g)$ is admissible.  A standard argument
\cite{pozniak,schwarz} then shows that $\partial^2=0$.

The homology of the Novikov complex depends only on the cohomology
class $[\alpha]$ and the covering $\tilde{X}$.  To describe it
topologically, choose a smooth triangulation of $X$, and lift the
simplices to obtain an equivariant triangulation of $\tilde{X}$.  We
denote the corresponding chain complex by $C_*(\tilde{X})$; this is a
module over the group ring $\Z[H]$.  There is then a natural
isomorphism
\begin{equation}
\label{eqn:novikov}
H_i(\Cnov_*)\simeq H_i(C_*(\tilde{X})\tensor_{\Z[H]}\Lambda).
\end{equation}
This was stated by Novikov \cite{novikov}, and proofs may be found in
\cite{pajitnov:old,pozniak,hl1}.  (Any cell decomposition would
suffice here, but we will shortly need to restrict to triangulations,
in order to define Reidemeister torsion refined by an Euler
structure.)

\subsection{Reidemeister torsion and the invariant $I$}
\label{sec:I}

The Novikov homology \eqref{eqn:novikov} often vanishes, at least
after tensoring with a field, and it is then interesting to consider
the Reidemeister torsion of the complexes $\Cnov_*$ and
$C_*(\tilde{X})$.

For certain rings $R$, including $\Z[H]$ and $\Lambda$, if $C_*$ is a
finite free chain complex over $R$ with a chosen basis $b$, then we
can define the {\bf Reidemeister torsion} $\tau(C_*)(b)\in Q(R)$, see
Appendix~\ref{app:torsion}.  The Novikov complex $\Cnov_*$ and
equivariant cell-chain complex $C_*(\tilde{X})$ have natural bases
consisting of lifts of the critical points or cells from $X$ to
$\tilde{X}$.  We denote the resulting torsion invariants by
\[
\overline{\Tmorse}\in Q(\Lambda)/\pm H,\quad\quad\quad
\overline{\Ttop}\in
Q(\Z[H])/\pm H.
\]

We have to mod out by $\pm H$ because of the ambiguity in choosing
lifts and ordering the bases.  It turns out that one can resolve the
former ambiguity by choosing an {\bf Euler structure}, see
Appendix~\ref{app:euler}.  The space $\Eul(X)$ of Euler structures is
a natural affine space over $H_1(X)$ defined by Turaev.  We thus
obtain refined torsion invariants, which are $H_1(X)$-equivariant maps
\begin{gather*}
\Tmorse: \Eul(X)\to Q(\Lambda)/\pm1,\\
\Ttop:\Eul(X)\to Q(\Z[H])/\pm1.
\end{gather*}

Results in \cite{turaev:euler} show that the refined topological
torsion $\Ttop$ depends only on the covering $\tilde{X}\to X$, and not
on the choice of triangulation.  For example, when $X$ is the
3-manifold obtained from 0-surgery on a knot $K$, the invariant
$\Ttop$ is related to the Alexander polynomial of $K$, see e.g.\
\cite{turaev:survey,hl1}.  By contrast, the Morse-theoretic torsion
$\Tmorse$ depends on the admissible pair $(\alpha,g)$, if $[\alpha]$
is fixed and nonzero. (See Example~\ref{ex:circle}.)  To get a
topological invariant, we must multiply by the zeta function.

\begin{definition}
\cite{hl2}
We define $\overline{I}\eqdef \overline{\Tmorse}\cdot\zeta\in
Q(\Lambda)/\pm H$, and
\[
I\eqdef \Tmorse\cdot\zeta:\Eul(X)\to Q(\Lambda)/\pm 1.
\]
\end{definition}

\begin{remark}
In the rest of this paper, in any equation involving the Reidemeister
torsions $\Ttop$ and $\Tmorse$ or the invariant $I$, it is to be
understood that there is an implicit `$\pm$' sign.  One can use a
homology orientation of $X$ to remove the sign ambiguity in the
topological torsion $\Ttop$ (see \cite{turaev:survey}), and presumably
in the Morse-theoretic torsion $\Tmorse$ as well, but we will not go
into that here.
\end{remark}

\subsection{The main results and basic examples}

Our main results are Theorems A and B below.  These were proved in
\cite{hl2} (generalizing \cite{hl1}) when the cohomology class of
$\alpha$ is rational.  A related result was proved by Pajitnov
\cite{pajitnov:new} at about the same time as \cite{hl2}; the
connection of this result with our work is discussed in
\S\ref{sec:conclusion}.

\bigskip

\noindent
{\bf Theorem A. } {\em For $(\alpha,g)$ admissible, the Morse theory
invariant $I$ depends only on the cohomology class $[\alpha]\in
H^1(X;\R)$ and the choice of covering $\tilde{X}$.  }

\bigskip

More specifically, the natural inclusion $\Z[H]\to\Lambda$ induces an
inclusion of quotient rings $\imath:Q(\Z[H])\longrightarrow
Q(\Lambda)$.  (To see this, one must check that the inclusion
$\Z[H]\to\Lambda$ sends nonzerodivisors to nonzerodivisors. This
follows from a ``leading coefficients'' argument or from
Lemma~\ref{lem:cyclotomic}.)  We then have:

\bigskip

\noindent
{\bf Theorem B. } {\em If $(\alpha,g)$ is admissible,
then the Morse theory invariant $I$ agrees with the topological
Reidemeister torsion:}
\[
\overline{I}=\imath(\overline{\Ttop})\in Q(\Lambda)/\pm H.
\]
We will also sketch a proof that the refined invariants agree,
i.e. that
\begin{equation}
\label{eqn:refined}
I=\imath\circ\Ttop:\Eul(X)\to Q(\Lambda)/\pm1.
\end{equation}
Of course \eqref{eqn:refined} implies Theorem A, since $\Ttop$ is a
topological invariant.  However our strategy will be to prove Theorem
A first, and then deduce Theorem B and \eqref{eqn:refined}.

\begin{example}
\label{ex:circle}
Suppose $X=S^1$ and $[\alpha]\neq 0$. We take $\tilde{X}=\R$, so that
$\Lambda\simeq\Z((t))$.  It is then not hard to check the following:
If $\alpha$ has no critical points, then $\overline{\Tmorse}=1$ and
$\zeta=(1-t)^{-1}$.  If $\alpha$ has critical points, then
$\overline{\Tmorse}=(1-t)^{-1}$ and $\zeta=1$.  In any case,
$\overline{\Ttop}=(1-t)^{-1}$.
\end{example}

\begin{example}
\label{ex:exact}
Suppose $\alpha=df$ with $f:X\to\R$ a Morse function. Then there are
no closed orbits, so $\zeta=1$. In this case it is classical that
$\overline{\Tmorse}=\overline{\Ttop}$, cf.\ \cite{milnor:whitehead}.
(Note that $\imath$ is the identity map in this case.)

Here is a sketch of a proof that in fact $\Tmorse=\Ttop$
(cf. \cite{hl2}).  Choose a triangulation $\mc{T}$ and let
$v_{\mc{T}}$ be the associated vector field as in
Appendix~\ref{app:euler}.  One can apparently find a Morse function
$f_{\mc{T}}$ and a metric $g_{\mc{T}}$ such that the gradient
$g_{\mc{T}}^{-1}df_{\mc{T}}$ is a perturbation of $v_{\mc{T}}$, so
that we have a natural isomorphism of chain complexes $\Cnov_*\simeq
C_*(\tilde{X})$, respecting the bases determined by an Euler
structure.  Then \eqref{eqn:refined} holds for
$(df_{\mc{T}},g_{\mc{T}})$, and by Theorem A it holds for all exact
1-forms.
\end{example}

\begin{example}
Suppose $\alpha=df$ where $f:X\to S^1$ is a fiber bundle with
connected fibers.  In particular, there are no critical points.  Let
$\tilde{X}$ be the infinite cyclic cover as in
Example~\ref{ex:cyclic}, so that $\Lambda\simeq\Z((t))$.  Let $\Sigma$
be a fiber, and let $\phi:\Sigma\to\Sigma$ be the return map obtained
by following the flow $-V$ from $\Sigma$ through $X$ and back to
$\Sigma$.  In this case the zeta function counts fixed points of
$\phi$ and its iterates with their Lefschetz signs:
\begin{equation}
\label{eqn:weil}
\zeta=\exp\sum_{k=1}^\infty\#\op{Fix}(\phi^k)\frac{t^k}{k}.
\end{equation}
There is a canonical Euler structure
$\xi_0=i_V^{-1}(0)$ (see Appendix~\ref{app:euler}), and
$\Tmorse(\xi_0)=1$.  One can also show (cf.\ \cite{fried:survey,hl2})
that
\[
\Ttop(\xi_0)=\prod_{i=0}^{n-1}\det(1-tH_i(\phi))^{(-1)^{i+1}}
\]
where $H_i(\phi)$ is the induced map on $H_i(\Sigma;\Q)$.  So Theorem
B gives here
\[
\exp\sum_{k=1}^\infty\#\op{Fix}(\phi^k)\frac{t^k}{k}
=\prod_{i=0}^{n-1}\det(1-tH_i(\phi))^{(-1)^{i+1}}. 
\]
This is equivalent to the Lefschetz theorem for $\phi^k$; to see this,
take the logarithmic derivative of both sides.  If we choose a larger
covering $\tilde{X}$, we obtain an equivariant version of the
Lefschetz theorem \cite{fried:survey,thesis}.
\end{example}

\begin{remark}
The relation between torsion and the zeta function in this example
goes back to Milnor \cite{milnor:cyclic}, and was extended to count
closed orbits of certain nonsingular flows by Fried
\cite{fried:homological}.  The version of the zeta function in
equation \eqref{eqn:weil} goes back to Weil \cite{weil}.
\end{remark}

\begin{example}
When $X$ is an oriented 3-manifold with $b_1(X)>0$, we conjectured in
\cite{hl1}, based on Taubes' work
\cite{taubes:sw=gr,taubes:counting,taubes:new}, that the invariant $I$
equals a certain reparametrization of the Seiberg-Witten invariant of $X$.  In
\cite{hl2}, we combined this conjecture with Theorem B to derive a
formula for the Seiberg-Witten invariant of $X$ in terms of
topological torsion, which had been conjectured by Turaev
\cite{turaev:spinc}.  This result was later independently proved by
Turaev \cite{turaev:new}, refining a result of Meng and Taubes
\cite{meng-taubes}, and indirectly verifying the conjecture in
\cite{hl1}.  For additional details and references see \cite{hl1,hl2}.

More recently, using ideas from TQFT, a paper by Donaldson
\cite{donaldson} has appeared giving an alternate approach to the
Meng-Taubes formula, and T.\ Mark \cite{mark} has given a more direct
proof of most of the conjecture in \cite{hl1}.
\end{example}

\subsection{Idea of the proofs}
\label{sec:outline}

The strategy for the proof of Theorem A is to analyze the effect on
$\Tmorse$ and $\zeta$ as we deform the pair $(\alpha,g)$, fixing the
cohomology class $[\alpha]$.  As long as the pair $(\alpha,g)$ remains
admissible, the Novikov complex and zeta function do not change.  But
in a generic 1-parameter family, the following types of bifurcations
(failures of admissibility) may occur:
\begin{description}
\item{(1)} A degenerate flow line from $\tilde{p}\in\Critt_i$ to
$\tilde{q}\in\Critt_{i-1}$,
\item{(2)} A degenerate closed orbit,
\item{(3)} A flow line from $\tilde{p}\in\Critt_i$ to
$\tilde{q}\in\Critt_i$, where $\pi(\tilde{p})\neq\pi(\tilde{q})$,
\item{(4)} A flow line from $\tilde{p}$ to $h\tilde{p}$, where $h\in
H$,
\item{(5)} Birth or death of two critical points at a degenerate
critical point.
\end{description}

The first two bifurcations change neither the Novikov complex nor the
zeta function.  This follows from compactness arguments for the moduli
spaces of closed orbits and flow lines.  Actually bifurcation (2)
includes two possibilities: simple cancellation of closed orbits of
opposite sign, and ``period-doubling bifurcations''.  Thus it is
important that the closed orbits are ``counted correctly'' by the zeta
function \eqref{eqn:zetaDef}, see \S\ref{sec:orbitCancelling}.

The third bifurcation does not affect the zeta function, but it does
change the Novikov complex, effectively performing a change of basis
in which $\tilde{p}$ is replaced by $\tilde{p}\pm\tilde{q}$.  This
does not change $\Tmorse$ because the change of basis matrix has
determinant one, cf.\ Proposition~\ref{prop:torsionAut}.

In the last two bifurcations, $\zeta$ and $\Tmorse$ both change, due
to the interaction of closed orbits with critical points.  In
bifurcation (4), a closed orbit homologous to $h$ is created or
destroyed, intuitively because the flow line from $\tilde{p}$ to
$h\tilde{p}$ is a ``broken closed orbit'', which should constitute a
boundary point in the one-dimensional moduli space of closed orbits as
time varies.  As a result, the zeta function is multiplied by a power
series $1\pm h+O(h^2)$.  At the same time, our understanding of
bifurcation (3) suggests that a change of basis occurs in the Novikov
complex in which $\tilde{p}$ is multiplied by a power series $1\pm
h+O(h^2)$.  Consequently the torsion $\Tmorse$ is multiplied by this
power series or its inverse.  We find in this way that $I$ is
unchanged ``to first order''.

The higher order terms are more difficult to understand, essentially
because the deformation upstairs in $\tilde{X}$ is not generic, due to
its $H$-equivariance, so that there are multiply broken closed orbits
and flow lines at the time of bifurcation.  It appears that $\zeta$
and $\Tmorse$ are simply multiplied by series of the form $(1\pm
h)^{\pm1}$.  But instead of trying to prove this, we consider a
non-equivariant perturbation of the deformation in a $k$-fold cyclic
cover $\hat{X}\to X$.  The idea is that invariance to first order in
$\hat{X}$ implies invariance to order $k$ in $X$, which we prove after
working out the behavior of the invariant $I$ with respect to finite
cyclic covers.  In this way we show that $I$ is unchanged by a
bifurcation of type (4).

Bifurcation (5) also has the subtlety of multiply broken flow lines
and closed orbits, arising from concatenations of flow lines from the
degenerate critical point to itself.  Here we use direct finite
dimensional analysis to show that every multiply broken closed orbit
or flow line leads to an unbroken closed orbit or flow line on the
side of the bifurcation time where the two critical points die, but
not on the other side.  Some miraculous algebra then yields that $I$
is invariant.

A small complication in the argument outlined above is that the times
at which bifurcations occur might not be isolated.  But bifurcations
involving ``short'' closed orbits or flow lines are isolated, and the
long bifurcations change only ``higher order'' terms in $I$. Taking a
limit in which we consider successively longer bifurcations, we
conclude that $I$ is invariant.

With Theorem A established, Theorem B can be deduced rather easily.
We already saw in Example~\ref{ex:exact} that Theorem B holds when
$\alpha$ is exact.  For general $\alpha$, we use a trick of F.\ Latour
which allows us to ``approximate'' $\alpha$ by an exact 1-form (!).
Namely, we let $f:X\to\R$ and replace $\alpha$ with the cohomologous
form
\[
\beta=\alpha+Cdf
\]
for $C\in\R$ large. The Novikov complex and zeta function of $\beta$
are the same as those of the rescaled form
\[
C^{-1}\beta = df+C^{-1}\alpha.
\]
We then see that for $C$ large, there are no closed orbits, and the
Novikov complex of $\beta$ coincides with that of $df$ (under an
inclusion of Novikov rings).  So by Example~\ref{ex:exact}, Theorem B
holds for $\beta$, and by invariance it holds for $\alpha$ as well.

\section{Proof of invariance I: preparation}
\label{sec:preparation}

In this section we make some general preparations for the proof of
Theorem A by bifurcation analysis.  In \S\ref{sec:bifurcationAnalysis}
we will undertake the analysis of specific bifurcations and complete
the proof of Theorem A.

\subsection{Semi-isolated bifurcations}
\label{sec:isolating}

Consider a 1-parameter family $\{(\alpha_t,g_t)\}$ of 1-forms and
metrics, parametrized by $t\in[0,1]$.  A generic family may have a
countably infinite set of bifurcations.  In this section we set up a
framework in which we only need to consider one bifurcation at a time.
More precisely, Lemma~\ref{lem:makeSense} makes sense of the change in
$I$ caused by a single bifurcation, and Lemma~\ref{lem:oneAtATime}
shows that if all these individual changes are zero, then $I$ is
invariant.  Note that we always assume that the cohomology class
$[\alpha_t]$ is fixed.

\begin{definition}
\label{def:broken}
A flow line between two critical points is {\bf degenerate} if it
corresponds to a nontransverse intersection of ascending and
descending manifolds.  A ($k$ times) {\bf broken flow line} from
$\tilde{p}$ to $\tilde{q}$ is a concatenation of flow lines from
$\tilde{p}$ to $\tilde{r_1}$ to $\tilde{r_2}$ to $\ldots$ to
$\tilde{r_k}$ to $\tilde{q}$, where $\tilde{r_1},\ldots,\tilde{r_k}$
are critical points and $k\ge 1$.  A {\bf broken closed orbit}
in the homology class $h$ is a (possibly broken) flow line from
$\tilde{p}$ to $h\tilde{p}$ for some critical point $\tilde{p}$.
\end{definition}

Let $\mc{M}_t(\tilde{p},\tilde{q})$ denote the space of (unbroken)
flow lines from $\tilde{p}$ to $\tilde{q}$ at time $t$.  Let
$\mc{O}_t(h)$ denote the space of (unbroken) closed orbits homologous
to $h$ at time $t$.  If the zeroes of $\alpha_t$ are nondegenerate for
all $t\in [t_1,t_2]$, then there is a canonical identification of
critical points $\Critt(t)=\Critt(t')$ for any $t,t'\in[t_1,t_2]$,
which we implicitly make below.

The following lemma implies that our invariant does not change if
there are no bifurcations, as a result of suitable compactness.

\begin{lemma}
\label{lem:unchanged}
Let $t_1<t_2$.  Suppose $(\alpha_t,g_t)$ is admissible at $t_1$ and $t_2$.
\begin{description}
\item{(a)}
Suppose there are no degenerate critical points for $t\in[t_1,t_2]$.
Let $\tilde{p}\in\Critt_i$ and $\tilde{q}\in\Critt_{i-1}$.  Suppose
there are no degenerate or broken flow lines from $\tilde{p}$ to
$\tilde{q}$ for $t\in[t_1,t_2]$.  Then
\[
\mc{M}_{t_1}(\tilde{p},\tilde{q})=\mc{M}_{t_2}(\tilde{p},\tilde{q}).
\]
\item{(b)}
Let $h\in H$.  Suppose there are no degenerate or broken closed orbits
homologous to $h$ for $t\in[t_1,t_2]$.  Then
\[
\mc{O}_{t_1}(h)=\mc{O}_{t_2}(h).
\]
\end{description}
Moreover, the above bijections are orientation preserving.
\end{lemma}

\begin{proof}
(a) It is enough to show that
$\bigcup_{t\in[t_1,t_2]}\mc{M}_t(\tilde{p},\tilde{q})$ is compact,
since it is then (possibly after a perturbation as in
\S\ref{sec:generic}) a one-manifold with boundary
$\mc{M}_{t_2}(\tilde{p},\tilde{q})-\mc{M}_{t_1}(\tilde{p},\tilde{q})$. Let
$S(\tilde{p})$ be a small sphere in the descending manifold around
$\tilde{p}$, and let $S(\tilde{q})$ be a small sphere in the ascending
manifold around $\tilde{q}$.  A flow line corresponds to a triple $(x,y,s)\in
S(\tilde{p})\times S(\tilde{q})\times\R$ such that downward flow from
$x$ for time $s$ hits $y$.  Compactness will follow from an upper
bound on $s$.  If $s$ is unbounded, then one can show as in
\S\ref{sec:zeta} that there is a broken or degenerate flow line from
$\tilde{p}$ to $\tilde{q}$ at some time $t\in[t_1,t_2]$.

For part (b), we get a similar compactness for
$\bigcup_{t\in[t_1,t_2]}\mc{O}_t(h)$.  Since the orbits remain
nondegenerate, none are created or destroyed, and the Lefschetz signs
cannot change.
\end{proof}

\begin{definition}
A {\bf bifurcation\/} of the family $\{(\alpha_t,g_t)\}$ is a time
$t_0\in\R$ such that the pair $(\alpha_{t_0},g_{t_0})$ fails to be
admissible.  The {\bf length} of a bifurcation $t_0$ is the smallest
of the following numbers:
\begin{description}
\item{(a)}
$0$, if $\alpha_{t_0}$ has a degenerate zero.
\item{(b)}
$[-\alpha_{t_0}](h)$, where $h$ is the homology class of a
degenerate or broken closed orbit.
\item{(c)}
$\int_\gamma-\alpha_{t_0}$, where $\gamma$ is a degenerate or
broken downward flow line.
\end{description}
\end{definition}

\begin{definition}
\label{def:good}
A time $t_0$ is {\bf good} if:
\begin{description}
\item{(a)}
$t_0$ is not a limit of bifurcations of bounded length.
\item{(b)}
For each $\epsilon>0$, the intervals $(t_0-\epsilon,t_0)$ and
$(t_0,t_0+\epsilon)$ both contain some times $t$ which are not bifurcations.
\end{description}
\end{definition}

\begin{definition}
\label{def:semiIsolated}
A bifurcation $t_0$ is {\bf semi-isolated} if $t_0$ is good, and:
\begin{description}
\item{(*)}
The pair $(\alpha_{t_0},g_{t_0})$ violates only one of the admissibility
conditions in Definition~\ref{def:admissible}, and in only one way.
\end{description}
\end{definition}

We now need to introduce the notion of limits in Novikov rings.
Given $x=\sum_ga_g\cdot g\in\Nov(G,N)$ and $R\in\R$, we write
``$x=O(R)$'' if $a_g=0$ whenever $N(g)<R$.  Given a sequence $\{x_n\}$ in
$\Nov(G;N)$ and $x\in\Nov(G;N)$, we write ``$\lim_{n\to\infty}x_n=x$''
if for every $R\in\R$ there exists $n_0$ such that $x-x_n=O(R)$ for
all $n>n_0$.

We can extend these definitions to the quotient ring $Q(\Lambda)$ as
follows.  If $G$ is a finitely generated abelian group, then by
Lemma~\ref{lem:cyclotomic} we have a decomposition
$Q(\Nov(G;N))=\oplus F_j$ into a sum of fields.  By
\eqref{eqn:sumsOfFields}, each field $F_j$ can be identified with the
tensor product of $\Nov(G/\Ker(N);N)$ with a certain field.  The
notion of ``$O(R)$'' is then well defined for elements of $F_j$. We
say that an element of $Q(\Nov(G;N))$ is ``$O(R)$'' if its projection
to each subfield $F_j$ is $O(R)$, and we define limits accordingly.

\begin{lemma}
\label{lem:makeSense}
If $t_0$ is good, then the limits as $t\nearrow t_0$ and $t\searrow
t_0$ of $\zeta$ and $(\Cnov_*,\partial)$ are well defined.
If moreover $t_0$ is not a bifurcation, then the left and right
limits of $\zeta$ and $\Cnov_*$ are equal to $\zeta(t_0)$
and $\Cnov_*(t_0)$.
\end{lemma}

\begin{proof}
Consider the limit as $t\nearrow t_0$.  There exists $\epsilon>0$ such
that all critical points of $\alpha_t$ are nondegenerate for
$t\in(t-\epsilon,t_0)$, so that $\Critt(t)=\Critt(t')$ for
$t,t'\in(t_0-\epsilon,t_0)$.  For convergence of $\partial$, we must
show that for $\tilde{p}\in\Critt_i$ and $\tilde{q}\in\Critt_{i-1}$,
there exists $x\in\Lambda$ such that
\[
\lim_{t\nearrow t_0}\sum_{h\in H}\langle\tilde{p},h\tilde{q}\rangle\cdot h=x,
\]
where $t$ ranges over any sequence of non-bifurcation values
converging to $t_0$ from below.  We use
Lemma~\ref{lem:unchanged}(a). For any path $\gamma$ from $\tilde{p}$
to $h\tilde{q}$, we have
\[
\int_\gamma-\alpha=C+[-\alpha](h)
\]
where $C$ is a constant which is independent of $h$ and varies
continuously with $t$.  Thus if $\gamma$ is a downward gradient flow
line and $[-\alpha](h)$ is bounded from above then
$\int_\gamma-\alpha$ is also bounded from above, so if we are
sufficiently close to $t_0$ then there are no degenerate or broken
flow lines from $\tilde{p}$ to $h\tilde{q}$ by definition of
semi-isolated, so $\langle\tilde{p},h\tilde{q}\rangle$ cannot change by
Lemma~\ref{lem:unchanged}(a).

Similary, Lemma~\ref{lem:unchanged}(b) implies that the zeta functions
converge.

The last sentence of the lemma follows from Lemma~\ref{lem:unchanged}.
\end{proof}

Let $\zeta^+$, $\zeta^-$, $(\Cnov_*^+,\partial^+)$,
$(\Cnov^-_*,\partial^-)$ denote these limits. An Euler structure gives
a basis for the limiting complexes $\Cnov_*^+$ and $\Cnov_*^-$; let
$\Tmorse^+$, $\Tmorse^-$ denote their Reidemeister torsion, and let
$I^\pm\eqdef\Tmorse^\pm\cdot\zeta^\pm$.

\begin{lemma}
\label{lem:limits}
If $t_0$ is good, then
\[
\lim_{t\nearrow t_0}I(t)=I^-(t_0),\;\;\;\;\;\;\;
\lim_{t\searrow t_0}I(t)=I^+(t_0),
\]
where $t$ ranges over non-bifurcations.
\end{lemma}

\begin{proof}
Consider the limit as $t\nearrow t_0$.  By definition we have
$\lim_{t\nearrow t_0}\zeta(t)=\zeta^-(t_0)$.  So we have to prove that
$\lim_{t\nearrow t_0}\Tmorse(t)(\xi)=\Tmorse^-(t_0)(\xi)$ where $\xi$
is a fixed Euler structure.

For $\epsilon$ sufficiently small we can identify the critical points
for different $t\in(t_0-\epsilon,t_0)$.  Fix a basis for $\Cnov_*$
consisting of a lift of each critical point to $\tilde{X}$, in the
equivalence class determined by $\xi$.

For a non-bifurcation $t$, recall that $\Tmorse$ is the sum of the
torsions of $\Cnov_*\tensor F_j$.  The torsion of $\Cnov_*\tensor F_j$
is zero if $\Cnov_*\tensor F_j$ is not acyclic; this criterion is
independent of $t$, by the Novikov isomorphism
\eqref{eqn:novikov}. Moreover, even if $t_0$ is a bifurcation, the
limiting complex $\Cnov_*^-\tensor F_j$ is acyclic if and only if
$\Cnov_*\tensor F_j$ is acyclic for all non-bifurcations $t$, because
the Novikov isomorphism \eqref{eqn:novikov}, as defined in \cite{hl1},
can be extended by a limiting argument to $\Cnov_*^-$.

When $\Cnov_*\tensor F_j$ is acyclic, we compute its torsion using
Proposition~\ref{prop:computeTorsion}.  Choose subbases $D_i$ and
$E_i$ as in Proposition~\ref{prop:computeTorsion} for
$\Cnov_*^-\tensor F_j$.  We can use these same subbases in the
interval $(t_0-\delta,t_0)$ for some $\delta$, because if the
determinants in Proposition~\ref{prop:computeTorsion} are nonzero in
the limiting complex, then they are nonzero near $t_0$.  The reason is
that each determinant for the limiting complex has a nonzero ``leading
term'' involving flow lines of length $<R$ for some $R$, which will be
unchanged near $t_0$ by Lemma~\ref{lem:unchanged}(a).

For $a,b\in F_j$ we have $\frac{1}{a}-\frac{1}{b}=O(R)$ when the
leading order of $a-b$ exceeds the leading order of $a$ and $b$ by at
least $R$.  This means that a high order change in the denominator of
$\Tmorse(\xi)$, as computed above, will change $\Tmorse(\xi)$ by high order
terms.  We are now done by condition \ref{def:good}(a) and
Lemma~\ref{lem:unchanged}(a).
\end{proof}

\begin{lemma}
\label{lem:oneAtATime}
Let $\{(\alpha_t,g_t)\}$ be a family parametrized by $t\in[0,1]$, with
$\alpha_t$ in a fixed cohomology class and $(\alpha_0,g_0)$ and
$(\alpha_1,g_1)$ admissible.  Suppose that every bifurcation
$t_0\in(0,1)$ is semi-isolated and satisfies $I^+(t_0)=I^-(t_0)$.
Then $I(0)=I(1)$.
\end{lemma}

\begin{proof}
Since every bifurcation satisfies condition \ref{def:good}(b), it
follows that the non-bifurcations are dense in $[0,1]$.  Moreover
every $t_0\in[0,1]$ is good.  (If $t_0$ is a non-bifurcation and fails
to satisfy condition \ref{def:good}(a), then a compactness argument
shows that $t_0$ is a bifurcation after all, giving a contradiction.)

By the assumptions and Lemma~\ref{lem:makeSense}, we have
$I^+(t_0)=I^-(t_0)$ for each $t_0\in[0,1]$.  It follows from
Lemma~\ref{lem:limits} that if we fix an Euler structure $\xi$ and
$R>0$, then for all $t_0\in[0,1]$, there exists $\epsilon>0$ such that
\[
I(t)(\xi)=I(t')(\xi)+O(R)
\]
for all non-bifurcations $t,t'\in(t_0-\epsilon,t_0+\epsilon)$.  Since
$[0,1]$ is compact and the non-bifurcations are dense, it follows that
$I(0)(\xi)-I(1)(\xi)=O(R)$. Taking $R\to\infty$, while keeping $\xi$
fixed, completes the proof.
\end{proof}

\subsection{Generic one-parameter families}
\label{sec:generic}

The following lemma implies that in a generic one-parameter family,
only the five types of bifurcations listed in \S\ref{sec:outline} may occur.

\begin{lemma}
\label{lem:genericDeformation}
Let $\{(\alpha_t,g_t),\;t\in[0,1]\}$ be a 1-parameter family with
$\alpha_t$ in a fixed cohomology class and $(\alpha_0,g_0)$ and
$(\alpha_1,g_1)$ admissible.  Then after a perturbation fixing the
endpoints, we may arrange that: 
\begin{description}
\item{(a)}
Near a degenerate critical point at time $t_0$, there are local
coordinates $x_1,\ldots,x_n$ in which 
\begin{equation}
\label{eqn:cerf}
g_t^{-1}\alpha_t=(x_1^2\pm (t-t_0),-x_2,\ldots,-x_i,x_{i+1},\ldots,x_n).
\end{equation}
\item{(b)}
Suppose that for $t\in[t_1,t_2]$, we have critical points
$\tilde{p}(t),\tilde{q}(t)\in\Critt(\alpha_t)$ which depend
continuously on $t$.  Then $\cup_t\mc{D}(\tilde{p})(t)$ and
$\cup_t\mc{A}(\tilde{q})(t)$ intersect transversely in
$X\times[t_1,t_2]$. 
\item{(c)} All bifurcations are semi-isolated.
\end{description}
If there are no degenerate critical points in the original family
$(\alpha_t,g_t)$, then we may choose this perturbation to be
$C^k$ small for any $k$. 
\end{lemma}

\begin{proof}
We will work with $C^k$ families.  After arranging (a), we will show
that in the space of $C^k$ families, there is a countable intersection
of open dense sets, whose elements are families with the desired
properties (b) and (c).  As in \cite{mcduff-salamon,taubes:counting},
we then obtain a dense set in $C^\infty$.

We begin by making the graph of $\cup_t\alpha_t$ transverse to
the 0-section of $T^*X\times[0,1]$.  Then $\cup_t\alpha_t^{-1}(0)$ is
a smooth 1-dimensional submanifold of $X\times[0,1]$.  We can further
arrange that $t$ is a Morse function on $\cup_t\alpha_t^{-1}(0)$ such
that all critical points have distinct values.  A critical point of
$t$ on $\cup_t\alpha_t^{-1}(0)$ is a pair $(x,t)$ where $\alpha_t$ has
a degenerate zero at $x$.  By a lemma of Cerf \cite{cerf} we can
choose (possibly time-dependent) local coordinates $x_1,\ldots,x_n$
near such a point so that
\[
\alpha_t=\frac{x_1^3}{3}\pm
(t-t_0)x_1-x_2^2-\cdots-x_i^2+x_{i+1}^2+\cdots+x_n^2.
\]
We now fix the metric $g_t$ on $X$ to be Euclidean near the origin in
these coordinates.  This gives (a). 

By a standard transversality argument, we can obtain (b) in a
countable intersection of open dense sets.  Fixing the metric near the
degenerate critical points does not interfere with the transversality
argument because no flow line or closed orbit is completely supported
near a degenerate critical point.

To obtain (c), we first arrange for the space of closed orbits to be
cut out transversely.  We then use a compactness argument to show that
(i) for each $R$, only finitely many bifurcations of length $<R$
occur.  We can arrange that these bifurcations occur at distinct
times, by intersecting with an open dense set of deformations.  So in
a countable intersection of open dense sets, we can arrange that (ii)
all bifurcations occur at distinct times. Now (i) and (ii) imply (c).
\end{proof}

\section{Proof of invariance II: bifurcation analysis}
\label{sec:bifurcationAnalysis}

In a generic one-parameter deformation given by
Lemma~\ref{lem:genericDeformation}, only the five types of
bifurcations listed in \S\ref{sec:outline} may occur, and all
bifurcations are semi-isolated.  In this section we will show that
$I^+=I^-$ for each bifurcation.  By Lemma~\ref{lem:oneAtATime}, this
will complete the proof of Theorem A.

\subsection{Cancellation of flow lines}

\begin{lemma}
\label{lem:flowLineCancelling}
Suppose $t_0$ is a semi-isolated bifurcation at which there is a
degenerate flow line from $\tilde{p}\in\Critt_i$ to
$\tilde{q}\in\Critt_{i-1}$. Then $\zeta^+(t_0)=\zeta^-(t_0)$ and
$(\Cnov_*^+,\partial^+)=(\Cnov_*^-,\partial^-)$.
\end{lemma}

\begin{proof}
By the definition of semi-isolated, we may choose $\epsilon>0$ such
that for all $t$ with $0<|t-t_0|<\epsilon$, there are no degenerate or
broken flow lines from $\tilde{p}$ to $\tilde{q}$. As in the proof of
Lemma~\ref{lem:unchanged}(a), the moduli space of flow lines from
$\tilde{p}$ to $\tilde{q}$ for $|t-t_0|<\epsilon$ is compact, so
$\langle\partial^-\tilde{p},\tilde{q}\rangle=
\langle\partial^+\tilde{p},\tilde{q}\rangle$, since the signed number
of boundary points of a compact 1-manifold is zero.

For every $R>0$, for every pair of critical points
$\tilde{r},\tilde{s}$ with index difference 1 and
$\int_{\tilde{s}}^{\tilde{r}}\alpha<R$, the coefficient
$\langle\partial\tilde{r},\tilde{s}\rangle$ likewise does not change
for $t$ sufficiently close to $t_0$.

For every $R>0$, the coefficients in the zeta function of $h$ with
$[-\alpha](h)<R$ do not change for $t$ sufficently close to $t_0$, by
Lemma~\ref{lem:unchanged}(b).
\end{proof}

\subsection{Cancellation of closed orbits}
\label{sec:orbitCancelling}

\begin{lemma}
Suppose $t_0$ is a semi-isolated bifurcation at which there
exists a degenerate closed orbit. Then
$(\Cnov_*^+,\partial^+)=(\Cnov_*^-,\partial^-)$ and
$\zeta^+=\zeta^-$.
\end{lemma}

\begin{proof}
The Novikov complex is unchanged as in the proof of
Lemma~\ref{lem:flowLineCancelling}.  To show that the zeta function
does not change, the idea is that locally the zeta function looks like
\eqref{eqn:weil}, and this is invariant because the signed number of
fixed points of a map is invariant, assuming suitable compactness.

More precisely, at time $t_0$ there is an isolated irreducible closed
orbit $\gamma$, with $[\gamma]=h$, such that $\gamma$ or some multiple
cover of it is degenerate.  Choose $x\in\gamma(S^1)\subset X$, and let
$D_\delta\subset X$ be a disc of radius $\delta$ transverse to
$\gamma$ and centered at $x$.  Let $\phi_{\delta,t}:D_\delta\to
D_\delta$ be the (partially defined) first return map for the flow
$g_t^{-1}\alpha_t$. We restrict the domain of $\phi_{\delta,t}$ to a
maximal connected neighborhood of $x$ on which it is continuous.
Define
\[
\zeta_{\delta,t}\eqdef
\exp\sum_{k=1}^\infty\frac{h^k}{k}\#\op{Fix}(\phi_{\delta,t}^k)
\]
for non-bifurcations $t$.  We claim that
\begin{equation}
\label{eqn:zetaDelta}
\frac{\zeta^+}{\zeta^-}= \lim_{\delta\to 0} \frac{\lim_{t\searrow
t_0}\zeta_{\delta,t}}{\lim_{t\nearrow t_0}\zeta_{\delta,t}}.
\end{equation}
To prove this, given $R>0$, we must find $\delta>0$ such that
$\frac{\zeta^+}{\zeta^-}=\frac{\lim_{t\searrow
t_0}\zeta_{\delta,t}}{\lim_{t\nearrow t_0}\zeta_{\delta,t}}+O(R)$.  By
the definition of semi-isolated, there exists $\epsilon>0$ so that for
$|t-t_0|<\epsilon$, all closed orbits $\gamma'$ with
$[-\alpha]([\gamma'])<R$ are nondegenerate, except for covers of
$\gamma$ at time $t=t_0$.  By compactness as in
Lemma~\ref{lem:unchanged}(b), we can choose $\delta$ sufficiently
small that no such closed orbit (other than covers of $\gamma$)
intersects $D_{2\delta}$ at time $t_0$.  Then for $|t-t_0|$
sufficiently small, the contribution to $\log\zeta$ from closed orbits
$\gamma'$ avoiding $D_\delta$ with $[-\alpha]([\gamma'])<R$ does not
change, and when moreover $t$ is not a bifurcation, the contribution
to $\log\zeta$ from all other closed orbits $\gamma'$ with
$[-\alpha]([\gamma'])<R$ is counted by the order $<R$ terms of
$\log\zeta_{\delta,t}$, as in \eqref{eqn:weil}.  This proves
\eqref{eqn:zetaDelta}.

Given any positive integer $k$, as above we can choose $\delta$ such
that at time $t_0$, no closed orbit $\gamma'$ with
$[-\alpha]([\gamma'])\le k[\-\alpha](h)$ (other than covers of
$\gamma$) intersects $D_{2\delta}$.  In particular, for $k'\le k$, the
boundary of the graph of $\phi_{\delta,t_0}^{k'}$ does not intersect
the diagonal in $D_\delta\times D_\delta$.  (Here we are compactifying
the graph as in the proof of Lemma~\ref{lem:powerSeries2}(b), see also
\cite{hl1}.)  It follows that $\#\op{Fix}(\phi_{\delta,t}^{k'})$ is
independent of $t$ for non-bifurcations $t$ close to $t_0$.  This
implies that
\[
\lim_{\delta\to 0} \frac{\lim_{t\searrow
t_0}\zeta_{\delta,t}}{\lim_{t\nearrow t_0}\zeta_{\delta,t}}=1.
\]
Together with \eqref{eqn:zetaDelta}, this proves the lemma.
\end{proof}

\begin{remark}
Here are two alternate approaches to proving this lemma, which might
generalize to Floer theory.

First, one might show that generically there is either a simple
cancellation of two orbits, or a ``period doubling'' bifurcation
corresponding to $(1+h)=(1-h^2)(1-h)^{-1}$ in the product formula
\eqref{eqn:productFormula}.  Related analysis appears in
\cite{taubes:counting} for the more complicated problem of counting
pseudoholomorphic tori in symplectic 4-manifolds.

Second, one might make the following heuristic rigorous.  For $h\in
H$, let $\mc{L}(h)$ denote the space of loops in $X$ homologous to
$h$, modulo reparametrization.  The coefficient of $h$ in $\log\zeta$,
\begin{equation}
\label{eqn:coefficient}
\sum_{\gamma\in\mc{O},\,[\gamma]=h}\frac{(-1)^{\mu(\gamma)}}{p(\gamma)}\in\Q,
\end{equation}
is formally the degree of a section of a vector bundle over $\mc{L}(h)$.
We divide by $p(\gamma)$ because $\mc{L}(h)$ is an orbifold with
$\Z/p$ symmetry around orbits with period $p$.  As long as there
is no interaction between closed orbits and critical points, so
that the zero set of the section remains compact, the coefficients
\eqref{eqn:coefficient}, and hence $\zeta$, should not change.
\end{remark}

\subsection{The slide bifurcation}

A {\bf slide bifurcation} is a semi-isolated bifurcation $t_0$ at
which there is a downward flow line from $\tilde{p}\in\Critt_i$ to
$\tilde{q}\in\Critt_i$.  (For real-valued Morse functions, this
bifurcation acts on the corresponding handle decomposition of $X$ by
sliding one handle over another.)  We assume that the flow line from
$\tilde{p}$ to $\tilde{q}$ is a transverse intersection of
$\cup_t\mc{D}(\tilde{p})(t)$ and $\cup_t\mc{A}(\tilde{q})(t)$.  

\begin{lemma}
\label{lem:slide}
For a slide bifurcation such that $\pi(\tilde{p})\neq\pi(\tilde{q})$
in $X$, we have
\begin{description}
\item{(a)}
$\zeta^+=\zeta^-$ and $\partial^+=A^{-1}\circ\partial^-\circ A$, where
$A:\Cnov_*\to\Cnov_*$ sends $\tilde{p}\mapsto\tilde{p}\pm\tilde{q}$
and fixes all other critical points $\tilde{s}$ with
$\pi(\tilde{s})\neq\pi(\tilde{p})$.
\item{(b)} In particular $I^+=I^-$.
\end{description}
\end{lemma}

\begin{proof}
For each flow line from $\tilde{s}\in\Critt_{i+1}$ to $\tilde{p}$ at
the bifurcation time, a flow line from $\tilde{s}$ to $\tilde{q}$ is
created or destroyed.  This follows from a standard gluing argument
\cite{floer} and can also be seen using finite dimensional methods as
in \cite{laudenbach}.  Similarly, for each flow line from $\tilde{q}$
to $\tilde{s}\in\Critt_{i-1}$ at the bifurcation time, a flow line
from $\tilde{p}$ to $\tilde{s}$ is created or destroyed.  By
Lemma~\ref{lem:unchanged}, no other flow lines or closed orbits are
created or destroyed or change sign.  This implies (a), after a check
that the orientations are consistent.  Part (b) follows by
Proposition~\ref{prop:torsionAut}, since $\det(A_i)=1$.
\end{proof}

\subsection{Torsion and zeta function of a finite cyclic cover}
\label{sec:cyclic}

We now digress to work out the behavior of the invariant $I$ with
respect to finite cyclic covers.  The answer is given in terms of the
Norm map from Galois theory.  This result will be needed when we use
nonequivariant perturbations in the next section, and may also be of
independent interest.

Suppose we have a short exact sequence of abelian groups
\[
0\longrightarrow K \stackrel{\imath}{\longrightarrow} H
\stackrel{m}{\longrightarrow} \Z/k\longrightarrow 0.
\]
Let $\rho:\hat{X}\to X$ be the $k$-fold cyclic covering whose
monodromy is the composition $\pi_1(X)\to H\stackrel{m}{\to}\Z/k$.
The covering $\tilde{X}\to X$ factors through $\rho$, and the covering
$\tilde{X}\to\hat{X}$ has automorphism group $K$.  We now want to
relate the invariants of $X$ and $\hat{X}$, choosing the covering
$\tilde{X}$ for both in Choice~\ref{choice:covering}.

We need the following algebraic notation.  Let
$\hat{\Lambda}\eqdef\op{Nov}(K;\imath^*[-\alpha])$.  The map $\imath$
induces a pushforward of Novikov rings
$\imath_*:\hat{\Lambda}\to\Lambda$ sending $\sum_{k\in K}a_k\cdot
k\mapsto \sum_{k\in K}a_k\cdot \imath(k)$.  Since $\imath$ has finite
kernel, there is also a pullback $\imath^*:\Lambda\to\hat{\Lambda}$
sending $\sum_{h\in H}a_h\cdot h\mapsto\sum_{k\in K}a_{\imath(k)}\cdot
k$.  The pushforward $\imath_*$ makes $\Lambda$ into a free module of
rank $k$ over $\hat{\Lambda}$.  If $y\in\Lambda$, then multiplication
by $y$ is an endomorphism of this module, whose determinant and trace
we denote by $\op{Norm}(y)$ and $\op{Tr}(y)$ respectively.

It will sometimes be convenient to assume that:
\begin{equation}
\label{eqn:assumption}
\mbox{$m$ annihilates the torsion subgroup of $H$.}
\end{equation}
In general, the map $\imath_*$ sends nonzerodivisors to
nonzerodivisors and hence induces a map on quotient rings
$Q(\hat{\Lambda})\to Q(\Lambda)$.  Recall from
Lemma~\ref{lem:cyclotomic}(a) that we have decompositions of
$Q(\hat{\Lambda})$ and $Q(\Lambda)$ into sums of fields.  Assumption
\eqref{eqn:assumption} implies that $\imath_*$ respects these
decompositions.  We then see from \eqref{eqn:sumsOfFields} that
$Q(\Lambda)$ is a free module of rank $k$ over $Q(\hat{\Lambda})$, so
$\op{Norm}$ extends to a multiplicative map $Q(\Lambda)\to
Q(\hat{\Lambda})$.

\begin{lemma}
\label{lem:NormTrace}
\begin{description}
\item{(a)}
If $y\in\Lambda$ then $\op{Tr}(y)=k\cdot\imath^*y$.
\item{(b)}
If $x\in\Lambda^+$, then
$
\;\log\op{Norm}(1+x)=\op{Tr}\log(1+x).
$
\item{(c)} Assuming \eqref{eqn:assumption}, if $y\in Q(\Lambda)$ and
$y\neq 0$, then $\op{Norm}(y)\neq 0$.
\end{description}
\end{lemma}

\begin{proof}
(a) is easy. To prove (b), let $\theta$ be a primitive $k^{th}$
root of $1$.  For $0\le i<k$, define a ring homomorphism
$\sigma_i:\Lambda\tensor\Z[\theta] \to \Lambda\tensor\Z[\theta]$ by
$\sigma_i(h\tensor 1)=h\tensor \theta^{i\cdot m(h)}$ for $h\in H$.  By
\cite[\S6.5]{lang}, we have
\begin{equation}
\label{eqn:NormTrace}
\op{Tr}(y) = \sum_{i=0}^{k-1}\sigma_i(y),\quad\quad
\op{Norm}(y) = \prod_{i=0}^{k-1}\sigma_i(y).
\end{equation}
The first of these identities implies that for $h\in H$,
\[
\op{Tr}(h)=\left\{\begin{array}{cl} kh & \mbox{if $m(h)=0$}\\
		0 & \mbox{if $m(h)\neq 0$}
		\end{array}\right.
\]
(which can also be seen more directly).  This proves (a).  To prove (b),
we compute
\[
\log\op{Norm}(1+x) =
\sum_{i=0}^{k-1}\log\sigma_i(1+x)
=\sum_{i=0}^{k-1}\sigma_i\log(1+x)
=\op{Tr}\log(1+x).
\]
Here the middle equality holds because $\log$ is defined using a power
series (see Remark~\ref{remark:powerSeries}) and $\sigma_i$ is a ring
homomorphism.  To prove (c), observe that assumption
\eqref{eqn:assumption} implies that $\sigma_i$ respects the field
decomposition of $Q(\Lambda)$.  Assertion (c) now follows from
\eqref{eqn:NormTrace} and the injectivity of $\sigma_i$.
\end{proof}

If $V$ is a vector field on $X$ with nondegenerate zeroes, then the
inverse image map $H_1(X,V)\to H_1(\hat{X},\rho^*V)$ induces a natural
pullback of Euler structures $\rho^*:\Eul(X)\to\Eul(\hat{X})$.  It is
clear that if $(\alpha,g)$ is admissible on $X$, then
$(\rho^*\alpha,\rho^*g)$ is admissible on $\hat{X}$.  We then have the
following result:

\begin{proposition}
\label{prop:cyclic}
\begin{description}
\item{(a)}
$\zeta(\hat{X})=\op{Norm}(\zeta(X))$.
\item{(b)} Under the assumption \eqref{eqn:assumption}, the following
diagram commutes:
\[
\begin{CD}
\Eul(\hat{X}) @>\Tmorse(\hat{X})>> Q(\hat{\Lambda})/\pm1\\
@A\rho^*AA @AA\op{Norm}A\\
\Eul(X) @>\Tmorse(X)>> Q(\Lambda)/\pm1.
\end{CD}
\]
\end{description}
\end{proposition}

\begin{proof}
(a) Every closed orbit $\hat{\gamma}$ in $\hat{X}$ is a lift of a
unique closed orbit $\gamma$ in $X$, with $[\gamma]\in K$.
Conversely, if $\gamma\in\mc{O}(X)$ and $[\gamma]\in K$, let
$\gamma_1$ denote the period one orbit underlying $\gamma$, and let $l$
be the order of $m([\gamma_1])$ in the group $\Z/k$.  Then $\gamma$
lifts to $k/l$ distinct closed orbits $\hat{\gamma}$, each of which
has period $p(\hat{\gamma})=p(\gamma)/l$ and Lefschetz sign
$(-1)^{\mu(\hat{\gamma})}=(-1)^{\mu(\gamma)}$.  Therefore
\[
\log\zeta(\hat{X}) =
\sum_{\hat{\gamma}\in\mc{O}(\hat{X})}
\frac{(-1)^{\mu(\hat{\gamma})}}{p(\hat{\gamma})}[\hat{\gamma}]
=\sum_{\gamma\in\mc{O}(X),\:[\gamma]\in K}
\frac{k(-1)^{\mu(\gamma)}}{p(\gamma)}[\gamma]
=k\imath^*\log\zeta(X).
\]
By Lemma~\ref{lem:NormTrace},
\[
k\imath^*\log\zeta(X)=\op{Tr}\log\zeta(X)=\log\op{Norm}\zeta(X).
\]
Combining the above equations and applying $\exp$ proves (a).

(b) A finite free complex $C_*$ over $\Lambda$ can be regarded as a
complex $\hat{C}_*$ over $\hat{\Lambda}$ with $k$ times as many
generators. Moreover, a basis $\{\lambda_1,\ldots,\lambda_k\}$ for
$\Lambda$ over $\hat{\Lambda}$ determines a map
$\phi:\Bas(C_*)\to\Bas(\hat{C}_*)$, and if $\chi(C_*)=0$ then the map
$\phi$ is independent of the choice of
$\{\lambda_1,\ldots,\lambda_k\}$.  Now we observe that if
$\xi\in\Eul(X)$, then the Novikov complex $\Cnov_*(\hat{X})$, with the
basis determined by $\rho^*\xi$, is obtained from $\Cnov_*(X)$ and
$\xi$ by this construction. So we need to show that
$\tau(\hat{C}_*)(\phi (b))=\op{Norm}(\tau(C_*)(b))$.  The assumption
\eqref{eqn:assumption} implies that $\imath_*$ and $\op{Norm}$ are
compatible with the decompositions of $Q(\hat{\Lambda})$ and
$Q(\Lambda)$ into sums of fields.  So we can restrict attention to a
complex $C_*\tensor F$ where $F\subset Q(\Lambda)$ is a field; let
$\hat{F}$ denote the corresponding field in $Q(\hat{\Lambda})$.  If
$C_*\tensor F$ is not acyclic, then $\hat{C}_*\tensor \hat{F}$ is not
acyclic either, so both torsions are zero.  If $C_*\tensor F$ is
acyclic, we can decompose it into a direct sum of 2-term acyclic
complexes.  Our claim then reduces to the fact that if $\partial$ is a
square matrix over $F$ and $\hat{\partial}$ is the corresponding
matrix over $\hat{F}$, then
$\det(\hat{\partial})=\op{Norm}(\det(\partial))$.  This follows from
the definition of $\op{Norm}$, after putting $\partial$ into Jordan
canonical form over an algebraic closure of $F$.
\end{proof}

\subsection{Sliding a critical point over itself}
\label{sec:bif4}

We now analyze bifurcation (4), in which a critical point slides over
itself, following the strategy described in \S\ref{sec:outline}.

If $p\in\Crit$ and $x\in\Lambda$, let $A_p(x):\Cnov_*\to\Cnov_*$
denote the $\Lambda$-module endomorphism which sends $\tilde{p}\mapsto
x\tilde{p}$ and fixes all other critical points $\tilde{s}$ with
$\pi(\tilde{s})\neq\pi(\tilde{p})$.

\begin{lemma}
\label{lem:powerSeries1}
Suppose $\tilde{p}\in\Critt_i$ slides over $h\tilde{p}$ for $h\in H$.  Then
\begin{description}
\item{(a)} There is a power series $x=1+\sum_{n=1}^\infty a_nh^n$, with
$a_n\in\Z$, such that
\begin{equation}
\label{eqn:Apx}
\partial^+=A_p(x)^{-1}\circ\partial^-\circ
A_p(x).
\end{equation}
\item{(b)}
In particular $\Tmorse^+=x^{(-1)^i}\cdot \Tmorse^-$.
\item{(c)}
The coefficient $a_1=\pm1$.
\end{description}
\end{lemma}

\begin{proof}
(a) Let $d$ denote the divisibility of $h$ in $H$.  (Note that $h$ is
not a torsion class.)  Let $k$ be a positive integer relatively prime
to $d$, and let $m:H\to \Z/k$ be a homomorphism sending $h\mapsto 1$.
Let $\rho:\hat{X}\to X$ be the $k$-fold cyclic cover with monodromy
$m$.  Then the critical points
$\tilde{p},h\tilde{p},\ldots,h^{k-1}\tilde{p}$ project to distinct
points in $\hat{X}$.

Let $R=[-\alpha](kh)$.  By semi-isolatedness, we can find $\epsilon>0$
such that no bifurcation of length $<R$ occurs between time
$t_0-\epsilon$ and $t_0+\epsilon$, other than the slide of $\tilde{p}$
over $h\tilde{p}$.  Choose a smaller $\epsilon$ if necessary so that
the pairs $(\alpha_{t_0\pm\epsilon},g_{t_0\pm\epsilon})$ are
admissible.  Perturb the pulled back family
$\{\rho^*(\alpha_t,g_t)|t\in[t_0-\epsilon,t_0+\epsilon]\}$, fixing the
endpoints, to satisfy the genericity conditions of
Lemma~\ref{lem:genericDeformation}.

By a compactness argument, we can choose the perturbation small enough
that no bifurcations of length $<R$ occur other than slides of
$h^i\tilde{p}$ over $h^j\tilde{p}$.  Then iterating Lemma~\ref{lem:slide}(a)
and using Lemma~\ref{lem:unchanged}(a), we find a power series
$x_k=1+\sum_{n=1}^{k-1}a_{n,k}h^n$ such that
\begin{equation}
\label{eqn:ank}
\partial^+=A(x_k)^{-1}\circ\partial^-\circ A(x_k) + O(R).
\end{equation}
(Here ``$O(R)$'' indicates a term involving flow lines $\gamma$ with
$\int_\gamma-\alpha\ge R$.)

Without loss of generality, $\partial^-\tilde{p}\neq 0$ or
$\langle\partial^-\tilde{s},\tilde{p}\rangle\neq 0$ for some $s$
(since otherwise equation \eqref{eqn:Apx} is vacuously true for any
$x$).  Then equation \eqref{eqn:ank} implies that for $n$ fixed,
$a_{n,k}$ is constant for large $k$.  If we define $a_n$ to be this
stable value of $a_{n,k}$, then equation \eqref{eqn:Apx} follows.

Assertion (b) follows from (a) and Proposition~\ref{prop:torsionAut}.

Now recall that the slide of $\tilde{p}$ over $h\tilde{p}$ comes from
a single transverse crossing of ascending and descending manifolds.
Under a sufficiently small perturbation of the deformation in
$\hat{X}$, this crossing will persist, and no other such crossing will
appear, by a compactness argument.  So for a sufficiently small
perturbation, $a_{1,k}=\pm 1$, and hence $a_1=\pm 1$.  This proves (c).
\end{proof}

\begin{lemma}
\label{lem:powerSeries2}
Suppose $\tilde{p}\in\Critt_i$ slides over $h\tilde{p}$.  Then
\begin{description}
\item{(a)} There is a power series
$y=1+\sum_{n=1}^\infty b_nh^n$ such that $\zeta^+=y\cdot\zeta^-.$
\item{(b)} $b_1=(-1)^{i+1}a_1$.
\end{description}
\end{lemma}

\begin{proof}
(a) By Lemma~\ref{lem:unchanged}(b), a closed orbit can be created or
destroyed in the bifurcation only if it is homologous to $kh$ for some
$k$.  So $\log(\zeta^+)-\log(\zeta^-)$ is a power series in $h$.  Thus
$\zeta^+/\zeta^-$ is a power series in $h$.  (A priori the
coefficients $b_n$ are rational; it's not important here, but we
actually know that $b_n\in\Z$, due to the product formula
\eqref{eqn:productFormula} for the zeta function.)

(b) Let $Z\subset X$ be a compact ``tubular'' neighborhood of the flow
line $\gamma$ from $p$ to itself at $t_0$.  There is a function
$f:Z\to \R/\Z$ such that $\alpha|_Z=\lambda df$ for some
$\lambda\in\R$.  Let $\Sigma\subset Z$ be a level set for $f$ away
from $p$.  The flow $-V$ induces a partially defined return map
$\phi:\Sigma\to\Sigma$.  Closed orbits homologous to $h$ near $\gamma$
are in one to one correspondence with fixed points of $\phi$.  A fixed
point of $\phi$ is an intersection of the diagonal
$\Delta\subset\Sigma\times\Sigma$ with the graph $\Gamma(\phi)$, and
the Lefschetz sign of the closed orbit equals the sign of the
intersection.  The graph $\Gamma(\phi)$ has a natural compactification
(see \cite{hl1}) to a manifold with corners $\overline{\Gamma}$ whose
codimension one stratum is
\[
\partial\overline{\Gamma}=(A(p)\times D(p))\cup Y.
\]
Here $D(p)$ and $A(p)$ are the ``first'' intersections of the
descending and ascending manifolds of $p$ with $\Sigma$, and $Y$ is a
component arising from trajectories that escape the neighborhood $Z$.
The number of closed orbits near $\gamma$ changes whenever $D(p)\times
A(p)$ crosses $\Delta$.  This is happening at time $t_0$ at a single
point, transversely, and an orientation check shows that the sign is
$(-1)^{i+1}a_1$.  No other closed orbits homologous to $h$ can be
created or destroyed, as in Lemma~\ref{lem:unchanged}(b).
\end{proof}

\begin{remark}
It should also be possible to prove (b) using a Floer-theoretic gluing
argument to show that in the homology class $h$, a single closed orbit
is created or destroyed.
\end{remark}

\begin{lemma}
\label{lem:bif4}
Suppose $\tilde{p}$ slides over $h\tilde{p}$.  Then
$I^+=I^-$.
\end{lemma}

\begin{proof}
By Lemmas~\ref{lem:powerSeries1} and \ref{lem:powerSeries2}, we can
write
\begin{equation}
\label{eqn:firstOrder}
I^+ = \left(\exp\sum_{n=2}^\infty
c_nh^n\right) I^-.
\end{equation}
for some $c_2,c_3,\ldots\in\Q$.  We need to show that each coefficient
$c_k$ vanishes.

Let $d$ denote the divisibility of $h$ in $H$.  Let $m:H\to\Z/dk$ be a
homomorphism which sends $h\mapsto d$ and annihilates the torsion
subgroup of $H$.  Let $\rho:\hat{X}\to X$ be the corresponding finite
cyclic cover.  By Proposition~\ref{prop:cyclic} and
Lemma~\ref{lem:NormTrace}(b),
\[
\begin{split}
I^+(\hat{X}) &= \op{Norm}\left(\exp\sum_{n=2}^\infty
c_nh^n\right)I^-(\hat{X})\\
&= \exp\left( dk\sum_{n=1}^\infty
c_{kn}h^{kn}\right) I^-(\hat{X}).
\end{split}
\]
As in Lemma~\ref{lem:limits}, we can choose $R$ sufficiently large
that a bifurcation of length $>R$ in $\hat{X}$ near $t_0$ will not
affect terms of order $[-\alpha](dkh)$ in $\Tmorse(\hat{X})$ or
$\zeta(\hat{X})$.  Now perturb the deformation in $\hat{X}$ as in the
proof of Lemma~\ref{lem:powerSeries1}, so that modulo bifurcations of
length $>R$, there are only slides of $h^i\tilde{p}$ over
$h^j\tilde{p}$.  When $k$ does not divide $j-i$, we know by
Lemma~\ref{lem:slide} that the torsion and zeta function in $\hat{X}$
do not change in such a slide.  When $j-i$ divides $k$, we apply the
analogue of \eqref{eqn:firstOrder} in the covering $\hat{X}$, to
conclude that $I(\hat{X})$ gets multiplied by $1+O(h^{2k})$.

It follows that $c_k=0$, as long as we know that $I^-(\hat{X})\neq 0$.  If
$\Cnov_*\tensor F$ is acyclic for at least one of the subfields $F$ of
$\Lambda$, then $I^{\pm}(X)\neq 0$, and it follows from
Lemma~\ref{lem:NormTrace}(c) and Proposition~\ref{prop:cyclic}(b) that
$I^{\pm}(\hat{X})\neq 0$, completing the proof.  If $\Cnov_*\tensor F$
is not acyclic for any $F$, then $I^{\pm}(X)=0$ and we have nothing to
prove.
\end{proof}

\begin{remark}
The last paragraph of the above proof could be avoided by working with
the {\em relative} torsion of the chain homotopy equivalence between
$\Cnov_*^-$ and $\Cnov_*^+$, cf.\ \S\ref{sec:conclusion}.
\end{remark}

\begin{remark}
A theorem of Shil'nikov \cite{arnold} asserts that in a generic
bifurcation of this type, a unique irreducible closed orbit is created
or destroyed.  By the product formula \eqref{eqn:productFormula},
$\zeta$ gets multiplied by $(1\pm h)^{\pm1}$.  By
Lemma~\ref{lem:bif4}, we see {\em a posteriori} that $\Tmorse$ is also
multiplied by such an expression.  A possible direct explanation for
this is that a flow line from $\tilde{p}$ to $h^n\tilde{p}$ is either
created for all $n$ or destroyed for all $n$.
\end{remark}

\subsection{Death of two critical points}
\label{sec:death}

We now analyze a semi-isolated {\bf death bifurcation} given by the
local model
\begin{equation}
\label{eqn:modelDeath}
V=(x_1^2+t-t_0,-x_2,\ldots,-x_i,x_{i+1},\ldots,x_n)
\end{equation}
in some neighborhood $U$ of the origin.  ({\bf Birth} is obtained from
death by reversing time.  Hence there is no loss of generality in
restricting attention to death.  However we will see below in
Proposition~\ref{prop:deathGeometry} that out of the death of two critical
points comes an abundance of new life.)

At time $t_0$ there is a
single degenerate critical point $r$. At time $t_0+\epsilon$, there
are no critical points in $U$. At time $t_0-\epsilon$, there are two
critical points $p=(-\sqrt{\epsilon},0,\ldots,0)$ and
$q=(\sqrt{\epsilon},0,\ldots,0)$ of indices $i$ and $i-1$
respectively.  Also there is a single downward gradient flow line in
$U$ from $p$ to $q$ in the positive $x_1$ direction, whose sign we
denote by $(-1)^\mu$.

If $x,y\in X$ are critical points of index difference one, let
$\mc{M}^-(x,y)$ denote the moduli space of flow lines from $x$ to $y$
immediately before the bifurcation.  If in addition $x,y$ are disjoint
from $p,q$, let $\mc{M}^+(x,y)$ denote the moduli space of flow lines
from $x$ to $y$ immediately after the bifurcation.  These moduli
spaces are well defined by the arguments in \S\ref{sec:isolating}.
Let $\mc{M}^0(r)$ denote the moduli space of flow lines from $r$ to
itself at the time of the bifurcation.  Let $\mc{O}^-$ and $\mc{O}^+$
denote the moduli spaces of closed orbits before and after the
bifurcation.

The following proposition says that for every (possibly multiply) broken
flow line or closed orbit at time $t_0$, a new flow line
or closed orbit is created after the two critical points die.

\begin{proposition}
\label{prop:deathGeometry}
\begin{description}
\item{(a)}
There is an orientation preserving bijection 
\[
\mc{O}^+=\mc{O}^-\bigcup\left(\bigcup_{k=1}^\infty (-1)^{\mu k+k+i+1}
\left(\mc{M}^0(r)\right)^{\times k}/(\Z/k)\right),
\]
which preserves total homology classes of orbits.  Here $\Z/k$ acts by
cyclic permutations.
\item{(b)}
If $x,y$ are critical points of index difference one which are
disjoint from $p,q$, then there is an orientation preserving bijection
\[
\begin{split}
\mc{M}^+(x,y)&=\mc{M}^-(x,y)\bigcup\\
	&\mc{M}^-(x,q)\times
		\bigcup_{k=0}^\infty(-1)^{(\mu+1)(k+1)}
		\left(\mc{M}^0(r)\right)^{\times k}
	\times\mc{M}^-(p,y)
\end{split}
\]
which preserves homology classes of flow lines.
\end{description}
\end{proposition}

\begin{proof}
In the calculations below, we will omit all orientations.

We first note that if $x,y$ are disjoint from $p,q$, then no flow
lines from $x$ to $y$ are destroyed, i.e. there is a natural inclusion
$\mc{M}^-(x,y)\to\mc{M}^+(x,y)$.  To see this, suppose to the contrary
that a flow line is destroyed.  Then by compactness there is a
sequence of flow lines from $x$ to $y$ before the bifurcation
converging to a broken or degenerate flow line from $x$ to $y$ at time
$t_0$.  There are no degenerate flow lines at $t_0$ (by the definition
of semi-isolated), so the limit flow line is broken, and the only
place it can be broken is at $r$.  In the neighborhood $U$, the broken
flow line approaches $r$ in the half space $(x_1>0)$ and leaves $r$ in
the half space $(x_1<0)$.  But such a broken flow cannot be the limit
as $\epsilon\to 0$ of unbroken flow lines at time $t_0-\epsilon$,
because there is a ``barrier'':  At time $t_0-\epsilon$, a downward
flow line cannot cross from $(x_1>\sqrt{\epsilon})$ to
$(x_1<\sqrt{\epsilon})$ within the neighborhood $U$, since the
downward gradient flow is in the positive $x_1$ direction for
$|x_1|<\sqrt{\epsilon}$.

Likewise, there is a natural inclusion $\mc{O}^-\to\mc{O}^+$.

To analyze what gets created, choose a small $\delta>0$ and let
$\Sigma_{\pm}\eqdef(x_1=\pm\delta)\subset U$.  Let
$D\eqdef\Sigma_-\cap\mc{D}(r)$ and $A\eqdef\Sigma_+\cap\mc{A}(r)$.
For $\epsilon$ small, let $f_\epsilon:\Sigma_+\to\Sigma_-$ denote the
partially defined map given by downward gradient flow at time
$t_0+\epsilon$.

Consider a broken closed orbit obtained by concatenating flow lines
$\gamma_1,\ldots,\gamma_k$ (in downward order) from $r$ to itself.
Choose $\delta$ small enough so that each $\gamma_i$ crosses $\Sigma_-$
immediately after leaving $r$ and crosses $\Sigma_+$ immediately
before returning.  Let $y_i\in D\subset\Sigma_-$ and $x_i\in
A\subset\Sigma_+$ denote the corresponding intersections of $\gamma_i$
with $\Sigma_{\pm}$.  The downward flow defines a return map $r_i$ from
a neighborhood of $y_i$ in $\Sigma_-$ to a neighborhood of $x_i$ in
$\Sigma_+$.

A new closed orbit approximating the broken one gets created for each
fixed point of the partially defined map
\begin{equation}
\label{eqn:partialMap1}
r_k\circ f_\epsilon\circ\cdots\circ r_1\circ f_\epsilon:
\Sigma_+\to\Sigma_+
\end{equation}
near $x_k$.  We will prove below that
\begin{equation}
\label{eqn:limitGraph1}
\lim_{\epsilon\to 0}
\Gamma(r_k\circ f_\epsilon\circ\cdots\circ r_1\circ f_\epsilon)
= A\times r_k(D).
\end{equation}
It follows that for $\epsilon$ small, the graph of
\eqref{eqn:partialMap1} intersects the diagonal once near $x_k\times
x_k$ transversely, because $A$ intersects $r_k(D)$ once transversely
at $x_k\times x_k$.  This proves (a).  (Note that no additional closed
orbits can be created, because by compactness a closed orbit can be
created only out of a broken closed orbit as above.)

To prove (b), suppose we have a broken flow line from $x$ to $y$ at
time $t_0$ consisting of a flow line $\gamma_0$ from $x$ to $r$,
followed by the concatenation of $\gamma_1,\ldots,\gamma_k$ and a
flow line $\gamma_{k+1}$ from $r$ to $y$.  Let $D'\subset \Sigma_+$
and $A'\subset\Sigma_-$ denote the corresponding intersections with
$\Sigma_+$ and $\Sigma_-$ of the descending manifold of $x$ and the
ascending manifold of $y$.  Let $\{x_0\}\eqdef\gamma_0\cap D'$ and
$\{y_{k+1}\}\eqdef\gamma_{k+1}\cap A'$.  A new flow line is created for each
intersection of the graph of the partially defined map
\begin{equation}
\label{eqn:partialMap2}
f_\epsilon\circ r_k\circ f_\epsilon\circ\cdots\circ r_1\circ
 f_\epsilon : \Sigma_+\to\Sigma_-
\end{equation}
with $D'\times A'$ near $x_0\times y_{k+1}$.  We will prove below that
\begin{equation}
\label{eqn:limitGraph2}
\lim_{\epsilon\to 0}
\Gamma(f_\epsilon\circ r_k\circ f_\epsilon\circ\cdots\circ r_1\circ
f_\epsilon) =
A\times D.
\end{equation}
It follows that for $\epsilon$ small, the graph of
\eqref{eqn:partialMap2} intersects $D'\times A'$ once transversely
near $x_0\times y_{k+1}$, because $A$ intersects $D'$ transversely at
$x_0$, and $D$ intersects $A'$ transversely at $y_{k+1}$.  This proves
(b).

We now prove equations \eqref{eqn:limitGraph1} and
\eqref{eqn:limitGraph2}.
We first note that by the local model \eqref{eqn:modelDeath}, we have
\begin{equation}
\label{eqn:limitGraph}
\lim_{\epsilon\to 0}\Gamma(f_\epsilon)=
A\times D\subset \Sigma_+\times\Sigma_-.
\end{equation}
In general, if $Y_1,Y_2,Y_3$ are manifolds and $\phi_1:Y_1\to Y_2$ and
$\phi_2:Y_2\to Y_3$ are any smooth maps, then $\Gamma(\phi_1)\times
Y_3$ intersects $Y_1\times\Gamma(\phi_2)$ transversely in $Y_1\times
Y_2\times Y_3$ and
\begin{equation}
\label{eqn:graphIntersection}
\Gamma(\phi_2\circ\phi_1)=\pi_{1,3}((\Gamma(\phi_1)\times Y_3)\cap
(Y_1\times\Gamma(\phi_2))),
\end{equation}
where $\pi_{1,3}:Y_1\times Y_2\times Y_3\to Y_1\times Y_3$ denotes the
projection.  Using \eqref{eqn:limitGraph} and
\eqref{eqn:graphIntersection} one proves \eqref{eqn:limitGraph1} and
\eqref{eqn:limitGraph2} together by induction on $k$.
\end{proof}

Let us now work out the algebraic consequences of the above lemma.
Choose lifts $\tilde{p}$ and $\tilde{q}$ of $p$ and $q$ which coalesce
at time $t_0$.  Choose a basis for $\Cnov_*^-$ so that $\tilde{p}$ and
$\tilde{q}$ are two of the basis elements.  For $\Cnov_*^+$, we can
use the same basis with $\tilde{p}$ and $\tilde{q}$ deleted.  Note
that these bases correspond to the same Euler structure, by
Definition~\ref{def:euler}.

In the former basis, we can write the
matrix for $\partial_i^-:\Cnov_i^-\to \Cnov_{i-1}^-$ in block form as
\begin{equation}
\label{eqn:blocks}
\partial_i^-=
\begin{pmatrix}
(-1)^\mu+\eta & v\\
w & N
\end{pmatrix}.
\end{equation}
Here $w$ is a column vector corresponding to $\tilde{p}$, and $v$ is a
row vector corresponding to $\tilde{q}$.  The power series $\eta$ counts
the flow lines in $\mc{M}^0(r)$ with their homology classes.  Note
that $\eta\in\Lambda^+$, so $(-1)^\mu+\eta$ is invertible.

We then have:

\begin{corollary}
\label{cor:deathAlgebra}
\begin{description}
\item{(a)}
$\Tmorse^-/\Tmorse^+=((-1)^\mu+\eta)^{(-1)^{i}}$.
\item{(b)}
$\zeta^+/\zeta^-=(1+(-1)^\mu\eta)^{(-1)^i}$.
\end{description}
\end{corollary}

\begin{proof}
By Proposition~\ref{prop:deathGeometry}(b), we have
$\partial_j^+=\partial_j^-$ for $j\neq i$, and
\[
\partial_i^+=N+\sum_{k=0}^\infty (-1)^{(\mu+1)(k+1)}w\eta^kv.
\]
We can rewrite this as
\begin{equation}
\label{eqn:partial_i}
\partial_i^+=N-w((-1)^\mu+\eta)^{-1}v.
\end{equation}

Now let $F$ be a subfield of $Q(\Lambda)$, as in
Lemma~\ref{lem:cyclotomic}.  Choose decompositions $\Cnov_*^+\otimes
F=D_*^+\oplus E_*^+$ as in Proposition~\ref{prop:computeTorsion}.  We
can then get subbases for $\Cnov_*^-\tensor F$ satisfying the
conditions of Proposition~\ref{prop:computeTorsion} by taking
$D_i^-=D_i^+\oplus\langle\tilde{p}\rangle$ and
$E_{i-1}^-=E_{i-1}^+\oplus\langle\tilde{q}\rangle$, and keeping the
other subbases fixed.  Let $N_s,v_s,w_s,\partial_s^\pm$ denote the
corresponding restrictions and/or projections of the $F$ components of
$N,v,w,\partial_i^\pm$.  Using \eqref{eqn:blocks} and
\eqref{eqn:partial_i}, we compute
\[
\begin{split}
\det(\partial_s^-:D_i^-\to E_{i-1}^-)
&=\det\begin{pmatrix} (-1)^\mu+\eta & v_s\\ w_s & N_s
\end{pmatrix}\\ & =
((-1)^\mu+\eta)\det(N_s-w_s((-1)^\mu+\eta)^{-1}v_s)\\
&= ((-1)^\mu+\eta)\det(\partial_s^+:D_i^+\to E_{i-1}^+).
\end{split}
\]
Putting this into Proposition~\ref{prop:computeTorsion} and summing
over subfields $F$, we obtain (a).  To prove (b), let us write
\[
\eta=\sum_{m=1}^\infty x_m\in\Lambda^+
\]
where there is one $x_m\in\pm H$ for each flow line from
$\tilde{r}$ to $h\tilde{r}$ at time $t_0$.  Then
\[
\begin{split}
\frac{\zeta^+}{\zeta^-} &= \exp\sum_{k=1}^\infty
\sum_{m_1,\ldots,m_k=1}^\infty
\frac{(-1)^{\mu k+k+i+1}}{k}x_{m_1}\cdots x_{m_k}\\
&= \left(1+(-1)^\mu\sum_{m}x_m\right)^{(-1)^i}.
\end{split}
\]
The first equality is a consequence of
Proposition~\ref{prop:deathGeometry}(a); the denominator $k$ arises because
summing over $k$-cycles and dividing by the period is equivalent to
summing over $k$-tuples and dividing by $k$.  The second equality can
be verified by taking the logarithm of both sides.  This proves (b).
\end{proof}

\begin{remark}
In the above calculation, we used the fact that
the determinant of a $2\times 2$ block matrix is given by
\[
\det\begin{pmatrix}\alpha&\beta\\\gamma&\delta\end{pmatrix}=
\det(\alpha)\det(\delta-\gamma\alpha^{-1}\beta),
\]
provided that $\alpha$ is invertible.  This identity played a key
role in \cite{hl2}, in a different argument.
\end{remark}

It follows from Corollary~\ref{cor:deathAlgebra} that $I$ is unchanged
under the death bifurcation, and this completes the proof of Theorem
A.

\section{Proof of Theorem B (comparison)}
\label{sec:comparison}

Let $(\alpha,g)$ be admissible.  We will now prove Theorem B,
identifying our invariant $I(\alpha,g)$ with topological
Reidemeister torsion.

We can reduce to the easier case of an exact one-form using the
following trick, which we learned from a
paper of Pajitnov
\cite{pajitnov:old}, who attributes it to F. Latour and J. Sikorav.
Choose $f:X\to\R$ such that $(df,g)$ is admissible, let $C\in\R$, and
define
\[
\beta\eqdef\alpha+Cdf.
\]

\begin{lemma}
\label{lem:latour}
If $C$ is sufficiently large, then $(\beta,g)$ is admissible, the
vector field $g^{-1}\beta$ has no nontrivial closed orbits, and there
is a canonical isomorphism of chain complexes
\begin{equation}
\label{eqn:beta}
\Cnov_*(\beta)= \Cnov_*(df)\otimes \Lambda
\end{equation}
respecting the bases determined by an Euler structure.
\end{lemma}

\begin{proof}
Since the Novikov complex is invariant under scaling, it makes no
difference if we take $\beta=df+\epsilon\alpha$ where $\epsilon$ is
small.   

Suppose $\gamma$ is a closed orbit of $g^{-1}\beta$.  The homology class
of $\gamma$ must be nonzero, since the cohomology class $[\alpha]$
pairs nontrivially with it.  We can then put a lower bound on the
length of $\gamma$ away from the critical points.  Since there is a
positive lower bound on $|df|$ away from the critical points, we
deduce a lower bound on $\int_\gamma(df+\epsilon\alpha)$.  If
$\epsilon$ is sufficiently small, then the closed orbit $\gamma$
cannot exist, or else we would get a positive lower bound on
$\int_\gamma df$, contradicting the fact that $\int_\gamma df=0$.

Transversality and intersection number are invariant under
small perturbations, so if $\epsilon$ is sufficiently small, then the
critical points of $\beta$ will be small perturbations of the critical
points of $f$ and remain nondegenerate, and the ascending and
descending manifolds will still intersect transversely with the same
intersection numbers.  This implies admissibility and \eqref{eqn:beta}.
\end{proof}

To prove Theorem B, choose a constant $C$ sufficiently large for the
conclusions of Lemma~\ref{lem:latour} to hold.  By Theorem A and
Lemma~\ref{lem:latour},
\begin{equation}
\label{eqn:B1}
I(\alpha,g)=I(\beta,g)=\Tmorse(\beta,g).
\end{equation}
We now use \eqref{eqn:beta} to relate $\Tmorse(\beta,g)$ to
$\Tmorse(df,g)$.  Note that the Novikov ring for $df$ is $\Z[H]$.  By
Lemma~\ref{lem:cyclotomic} we have decompositions
\[
\begin{split}
Q(\Z[H]) &= \bigoplus_{j=1}^mF_j,\\
Q(\Lambda) &=\bigoplus_{j=1}^mF_j'
\end{split}
\]
into sums of fields such that $\imath(F_j)\subset F_j'$, where
$\imath:Q(\Z[H])\to Q(\Lambda)$ is the natural inclusion.  By
Proposition~\ref{prop:computeTorsion} we see that $\Cnov_*(df)\tensor F_j$ is
acyclic if and only if $\Cnov_*(df)\tensor F_j'$ is, and by
\eqref{eqn:beta},
\begin{equation}
\label{eqn:B2}
\Tmorse(\beta,g)=\imath\circ \Tmorse(df,g):\Eul(X)\to
\frac{Q(\Lambda)}{\pm1}. 
\end{equation}
By Example~\ref{ex:exact},
\begin{equation}
\label{eqn:B3}
\overline{\Tmorse}(df,g)=\overline{\Ttop}.
\end{equation}
Equations \eqref{eqn:B1}, \eqref{eqn:B2}, and \eqref{eqn:B3} prove
Theorem B.

\begin{remark}
D. Salamon points out that instead of Lemma~\ref{lem:latour}, one can
use a lemma of Pozniak \cite{pozniak} asserting that for any
cohomology class $a\in H^1(X;\R)$, there are admissible pairs
$(g_1,\alpha)$ and $(g_2,df)$, where $[\alpha]=a$, with identical
vector fields $g_1^{-1}\alpha=g_2^{-1}df$.
\end{remark}

\begin{remark}
A rigorous justification of the sketch in Example~\ref{ex:exact} would
allow us to remove the bars from \eqref{eqn:B3} and deduce the
refinement by Euler structures \eqref{eqn:refined}.
\end{remark}

\section{Conclusion}
\label{sec:conclusion}

There are several directions in which the results of this paper might
be generalized.

\paragraph{Algebraic refinements.}
 There are sharper notions of torsion which are defined less often.
The sharpest is Whitehead torsion \cite{cohen,milnor:whitehead}, which
is only defined for an acyclic complex over a ring $R$, and lives in the ring
$\overline{K_1}(R)$.  One can also define the ``relative'' Whitehead
torsion of a chain homotopy equivalence between two complexes which
need not be acyclic.  A homotopy $\{(\alpha_t,g_t)\}$ between
admissible pairs $(\alpha_0,g_0)$ and $(\alpha_1,g_1)$, with the
cohomology class $[\alpha_t]$ fixed, induces a chain homotopy
equivalence between the two Novikov complexes, via ``continuation''
(cf. \cite{pozniak,schwarz}).  It should be possible to upgrade the
algebra in Theorem A to show that {\em the Whitehead torsion of the
continuation map equals the ratio of the two zeta functions}.  Modulo
Euler structures, and under slightly stronger genericity assumptions,
this follows {\em a posteriori} from the paper of Pajitnov
\cite{pajitnov:new}.

One might also generalize our results to nonabelian covers.  We
believe that if such a generalization exists, then the bifurcation
analysis in this paper should suffice to prove it.  The difficulty
seems to be to formulate a result.  In this direction, several earlier
works, including \cite{sikorav,latour,pajitnov:old,pajitnov:surgery},
investigated the Novikov complex for the universal cover and its
Whitehead torsion; zeta functions for the universal cover were
introduced in \cite{nicas}.

\paragraph{Infinite dimensions.}
Floer theory considers finite dimensional moduli spaces of flow lines
of closed 1-forms on certain infinite dimensional manifolds.  Several
people have suggested to us that for any such setup, one can at least
formally define an analogue of our invariant $I$.  Theorem A might
generalize to prove that such a construct is invariant under exact
deformations.  (Whitehead torsion in Floer theory, without the zeta
function, is studied in \cite{fukaya,sullivan}.)

To give one example, consider the Floer theory of a symplectomorphism
$f:X\to X$ of a symplectic manifold $X$.  Let
$M_f\eqdef X\times[0,1]/(x,0)\sim(f(x),1)$ denote the mapping torus.  One
defines a complex $CF_*(X,f)$ whose chains are fixed points of $f$ and
whose boundary operator counts pseudoholomorphic annuli in
$M_f\times\R$ which converge at either end to loops coming from fixed
points.  One can define the algebraic Reidemeister torsion of this
complex just as in the finite dimensional case.  Furthermore the
analogue of the zeta function should count certain pseuodholomorphic
tori in $M_f\times S^1(r)$, where $S^1(r)$ is the circle of radius
$r$.  The signs of the tori can be defined using spectral flow,
cf. \cite{taubes:counting}.  Due to the $S^1$ action, to get a moduli
space of expected dimension zero, we must allow $r$ to
vary.  During a deformation, tori may disappear if $r\to\infty$.
However the energy of a long torus will be small on most of it, so
part of the torus should be approaching a critical point, in which
case we expect the loss of the torus to be reflected in a change in
torsion as in bifurcations (4) and (5) on the list in
\S\ref{sec:outline}.

We have tried to write the proof of Theorem A in such a way that it
can be easily generalized to Floer theory.  However a better
understanding is needed of the gluing of multiply broken flow lines,
which arises in bifurcations (4) and (5).  In particular one would
like to understand: On what {\em side} of the bifurcation time are
things created or destroyed?  The ``nonequivariant perturbation''
trick, which we used to evade this issue in bifurcation (4), does not
appear to work for bifurcation (5), where we resorted in this paper to
purely finite-dimensional methods.

We remark that Floer proved invariance of Floer homology by directly
constructing a chain homotopy equivalence, without using bifurcation
analysis.  It seems however that bifurcation analysis is necessary to
prove the invariance of torsion; roughly, one needs to see that the
chain homotopy equivalence is composed out of a restricted set of matrix
operations.

\paragraph{Other vector fields.}

The fact that our vector field $V$ is dual to a closed 1-form is used
mainly to give uniform bounds on the numbers of closed orbits and flow
lines so that finite counting is possible.  Fried
\cite{fried:homological} relates zeta functions to Reidemeister
torsion for a rather different kind of vector field, assuming that
there are no critical points.  We do not know to what class of vector
fields our results can be generalized.  In the setting of
combinatorial Morse theory, a statement resembling Theorem B was
recently proved by Forman \cite{forman}.

\paragraph{Acknowledgments.}
This paper is a revised and abridged version of my PhD thesis
\cite{thesis}, which was supported by a Sloan Dissertation Fellowship,
and earlier by a National Science Foundation Graduate Fellowship.  I
am extremely grateful to my advisor Cliff Taubes for his support and for
the excellent suggestion which got this project started.  I am
indebted to Yi-Jen Lee for teaching me many things during our
collaboration on \cite{hl1,hl2}, and for pointing out some mistakes in
\cite{thesis}. I thank R. Bott, R. Forman, D. Fried, D. Salamon,
M. Schwarz, P. Seidel, and J. Weber for helpful conversations and
encouragement, and K. Conrad for help with Galois theory.  I thank
A. Pajitnov and V. Turaev for sending me their highly relevant papers.

\begin{appendix}

\section{The algebra of Reidemeister torsion}
\label{app:torsion}

In this appendix we review the algebra that underlies the definitions of
topological and Morse-theoretic Reidemeister torsion, and which is
needed starting in \S\ref{sec:I}.

We call a complex $(C_i,\partial)$ over a ring $R$ {\bf free} if each
$C_i$ is a free $R$-module, and {\bf finite} if
$\sum_i\op{rk}(C_i)<\infty$.  A {\bf basis} $b$ of a finite free
complex consists of an ordered basis $b_i$ for each $C_i$.  We declare
two bases $b,b'$ to be equivalent if
\[
\prod_{i}[b_{2i},b'_{2i}]=\prod_i[b_{2i+1},b'_{2i+1}]\in R,
\]
where $[b_i,b_i']\in R$ denotes the determinant of the change of basis
matrix from $b_i$ to $b_i'$.  (We assume that the bases $b_i$ and
$b_i'$ have the same cardinality, which could fail for pathological
rings $R$.)  We denote the set of equivalence classes by $\Bas(C_*)$.

If $(C_*,\partial)$ is a finite complex over a field $F$, we define the
Reidemeister torsion
\[
\tau(C_*,\partial):\Bas(C_*)\to F
\]
as follows.  The standard short exact sequences $0\to
Z_i\to C_i\stackrel{\partial}{\to} B_{i-1}\to 0$ and $0\to B_i\to
Z_i\to H_i\to 0$ give rise to isomorphisms
\[
\begin{split}
\det(C_i)&\too\det(Z_i)\otimes\det(B_{i-1}),\\
\det(Z_i)&\too\det(B_i)\otimes\det(H_i),
\end{split}
\]
where `$\det$' denotes top exterior power.  Putting the second
isomorphism into the first gives an isomorphism
\[
\det(C_i)\too\det(H_i)\otimes\det(B_i)\otimes\det(B_{i-1}).
\]
When we take the alternating product over $i$, the $B$'s cancel and we
obtain an isomorphism
\begin{equation}
\label{eqn:determinantIsomorphism}
\Bas(C_*)=\bigotimes_i\det(C_i)^{\tensor(-1)^i}
\too\bigotimes_i\det(H_i)^{\tensor(-1)^i}.
\end{equation}

\begin{definition}
If $(C_*,\partial)$ is acyclic, then
$\bigotimes_i\det(H_i)^{\tensor(-1)^i}=F$, and we define the {\bf
Reidemeister torsion} $\tau(C_*,\partial)$ to be the map
\eqref{eqn:determinantIsomorphism}.  If $(C_*,\partial)$ is not
acyclic, we define $\tau(C_*,\partial)\eqdef 0$.
\end{definition}

In practice, one can compute torsion as an alternating product of
determinants of square submatrices of $\partial$.  More precisely:

\begin{proposition}
\label{prop:computeTorsion}
Let $(C_i,\partial)$ be a finite acyclic complex over a field $F$ with
a fixed basis $b$.  We can find decompositions $C_i=D_i\oplus E_i$
such that:
\begin{description}
\item{(i)}
$D_i$ and $E_i$ are spanned by subbases of $b_i$, and
\item{(ii)}
The map $\partial_s\eqdef\pi_{E_{i-1}}\circ\partial|_{D_i}:D_i\to E_{i-1}$ is
an isomorphism.
\end{description}
We then have
\[
\tau(C_*,\partial)(b) =
\pm\prod_{i}\det(\partial_s:D_i\to E_{i-1})^{(-1)^i}
\]
where the determinants are computed using the subbases of $b$.
\noProof
\end{proposition}

Suppose now that $C_*$ is a finite free complex over a ring $R$, such
that the total quotient ring $Q(R)$ is a finite direct sum of fields,
\begin{equation}
\label{eqn:sumOfFields}
Q(R)=\bigoplus_jF_j.
\end{equation}

\begin{definition}
\cite{turaev:spinc}
Under the above assumption, we define
\begin{align*}
\tau(C_*,\partial):\Bas(C_*)&\too Q(R),\\
b &\longmapsto
\sum_j\tau(C_*\otimes_RF_j,\partial\tensor1)(b\otimes 1).
\end{align*}
This depends only on $R$, i.e. the decomposition
\eqref{eqn:sumOfFields} is unique, because the fields $F_j$ are
characterized as the minimal ideals in $Q(R)$.
\end{definition}

This definition applies to the complexes of interest in this paper, by:

\begin{lemma}
\label{lem:cyclotomic}
Let $G$ be a finitely generated abelian group.  Then:
\begin{description}
\item{(a)}
The total quotient rings of $\Z[G]$ and $\Nov(G;N)$ are finite sums of fields.
\item{(b)}
These decompositions are compatible with the inclusion $\Z[G]\to\Nov(G;N)$.
\end{description}
\end{lemma}

\begin{proof}
(cf. \cite{turaev:spinc}) Choose a splitting $G=K\oplus F$ where $K$
is finite and $F$ is free.  Then $\Z[G]=\Z[K]\otimes\Z[F]$ and
$\Nov(G;N)=\Z[K]\otimes \Nov(F;N)$.  The total quotient ring of
$\Z[K]$ is a finite sum of (cyclotomic) fields,
$
Q(\Z[K])=\bigoplus_jL_j.
$
We then have
\begin{equation}
\label{eqn:sumsOfFields}
\begin{split}
Q(\Z[G])&=\bigoplus_jL_j\tensor Q(\Z[F]),
\\
Q(\Nov(G;N))&=\bigoplus_jL_j\tensor Q(\Nov(F;N)).
\end{split}
\end{equation}
A ``leading coefficients'' argument shows that $\Z[F]$ and $\Nov(F;N)$
are integral domains, so $L_j\tensor Q(\Z[F])$ and $L_j\tensor
Q(\Nov(F;N))$ are fields.  Thus equation \eqref{eqn:sumsOfFields}
proves (a) and (b).
\end{proof}

The following ``change of basis''  formula is important in
\S\ref{sec:bifurcationAnalysis}.

\begin{proposition}
\label{prop:torsionAut}
Let $(C_*,\partial)$ be a finite free complex over $R$, where $Q(R)$
is a finite sum of fields.  If $A_*\in\op{Aut}(C_*)$ preserves the
grading, then
\begin{equation}
\tau(C_*,A^{-1}\partial A)
=\tau(C_*,\partial)\cdot\prod_i\det(A_i)^{(-1)^i}.
\tag*{$\Box$}
\end{equation}
\end{proposition}

\section{Euler structures}
\label{app:euler}

In this appendix we explain how to resolve the $H$ ambiguity in
topological and Morse-theoretic Reidemeister torsion (cf.\
\S\ref{sec:I}), using Turaev's Euler structures.

We begin with a definition of Euler structures which is slightly
different from Turaev's.  If $v$ is a smooth vector field on $X$ with
nondegenerate zeroes, let $\Eul(X,v)$ denote the set of homology
classes of 1-chains $\gamma$ with $\partial\gamma=v^{-1}(0)$, where
$v^{-1}(0)$ is oriented in the standard way.  The set $\Eul(X,v)$ is a
subset of the relative homology $H_1(X,v^{-1}(0))$, and it is an
affine space modelled on $H_1(X)$.  The set $\Eul(X,v)$ is nonempty
because we are assuming $\chi(X)=0$.

If $v_0,v_1$ are two such vector fields, define
\[
\phi_{v_1,v_0}:\Eul(X,v_0)\to\Eul(X,v_1)
\]
as follows.  Let $w$ be a vector field on $X\times [0,1]$ such that
$v_i=w|_{X\times\{i\}}$ and $w^{-1}(0)$ is cut out transversely.  The
orientation convention gives $\partial
w^{-1}(0)=v_1^{-1}(0)-v_0^{-1}(0)$.  Suppose $\gamma\in\Eul(X,v_0)$.
Since $H_1(X\times[0,1],X\times\{1\})=0$, there is a 2-chain
$\Sigma\subset X\times [0,1]$ with
$\partial\Sigma=-w^{-1}(0)-\gamma\;(\op{rel} X\times\{1\})$.  We
define $\phi_{v_1,v_0}(\gamma)\eqdef \partial\Sigma+w^{-1}(0)+\gamma$.

\begin{definition}
\label{def:euler}
One can check that (a) $\phi_{v_1,v_0}$ is independent of $w$ and
$\Sigma$, (b) $\phi_{v,v}=\op{id}$, and (c)
$\phi_{v_2,v_0}=\phi_{v_2,v_1}\phi_{v_1,v_0}$.  This implies that all
the spaces $\Eul(X,v)$ are canonically isomorphic to a single affine
space over $H_1(X)$.  We denote this space by $\Eul(X)$ and call an
element of it an {\bf Euler structure}.  We let $i_v:\Eul(X)\to\Eul(X,v)$
denote the canonical isomorphism.
\end{definition}

It should be emphasized that the affine space $\Eul(X)$ is not
canonically isomorphic to $H_1(X)$. For example, when $v_0,v_1$ have
no zeroes, the map $\phi_{v_1,v_0}$ does not necessarily respect the
identifications $\Eul(X,v_i)\simeq H_1(X)$.

\begin{remark}
When $\dim(X)>1$, Turaev \cite{turaev:euler} defines a (smooth) Euler
structure to be a nonsingular continuous vector field, modulo homotopy
through vector fields which remain nonsingular in the complement of a
ball during the homotopy.  To go from our definition to Turaev's,
represent $\gamma\in\Eul(X,v)$ by disjoint paths connecting the zeroes
of $v$, and cancel the zeroes of $v$ in a neighborhood of $\gamma$.
\end{remark}

We now explain how Euler structures determine (equivalence classes of)
bases for the Novikov complex.

\begin{definition}
We define a map
\begin{equation}
\label{eqn:bmorse}
\Eul(X)\too\Bas(\Cnov_*)/\pm1
\end{equation}
as follows.  If there are no
critical points, then $\Cnov_i=\{0\}$, so $\Bas(\Cnov_*)=H_1(X)$.  In this
case we define the map \eqref{eqn:bmorse} to be the composition
$\Eul(X)\stackrel{i_V}{\to}\Eul(X,V)=H_1(X)$.

If $V^{-1}(0)\neq\emptyset$, then given $\xi\in\Eul(X)$, we can
represent $i_V(\xi)\in \Eul(X,V)$ by a chain $\gamma$ consisting only
of paths connecting the zeroes of $V$, such that each critical point
is in one component of $\gamma$. Choose a lift $\tilde{\gamma}$ of
$\gamma$ to $\tilde{X}$.  The induced lifts of the zeroes of $V$ to
the endpoints of $\tilde{\gamma}$ determine a basis for $\Cnov_*$.
\end{definition}

The equivalence class of this basis does not depend on the choice of
lift $\tilde{\gamma}$, because the boundary of each component of
$\gamma$ consists of two critical points whose indices have opposite
sign.  It is also independent of $\gamma$.

We now consider bases of the equivariant cell complex, along the lines
of \cite{turaev:euler}.  There is a standard vector field $v_i$ on the
standard $i$-simplex with a sink at the center of the simplex, with no
other zeroes in the interior, which restricts to $v_j$ on each
$j$-dimensional face, and which points inward near the boundary
\cite{turaev:euler}.  Putting the vector fields $v_i$ onto the
simplices of our triangulation $\mc{T}$, we obtain a continuous vector
field $v_{\mc{T}}$ on $X$.  We can perturb this to a smooth vector
field $v$ with a nondegenerate zero of sign $(-1)^i$ in the center of
each $i$-simplex.

\begin{definition}
We define a map
\[
\Eul(X)\too\Bas(C_*(\tilde{X}))/\pm1
\]
as follows.  Given $\xi\in\Eul(X)$, represent $i_v(\xi)\in
\Eul(X,v)$ by a chain $\gamma$ consisting only of paths connecting the
centers of the simplices in pairs.  Choose a lift $\tilde{\gamma}$ of
$\gamma$ to $\tilde{X}$. Each simplex $\sigma$ in $X$ now has a unique
lift in $\tilde{X}$ such that the center of $\sigma$ is lifted to one
of the points of $\partial\tilde{\gamma}$. These simplices in
$\tilde{X}$ give a basis for $C_*(\tilde{X})$.
\end{definition}

The equivalence class of this basis does not depend on the
perturbation $v$, the path $\gamma$, or the lift $\tilde{\gamma}$.

\end{appendix}

\end{document}